\documentclass[12pt,reqno]{amsart}
\usepackage{amsmath,amssymb,mathrsfs,amsfonts,mathabx}
\usepackage{graphicx,cite,cases}
\usepackage{xcolor}
\usepackage{epstopdf}

\newtheorem{theorem}{Theorem}[section]

\numberwithin{equation}{section}

\begin{document}

\title[an inverse source problem]{An inverse random source problem
for the Helium production-diffusion equation driven by a fractional Brownian motion}

\author{Jing Li}

\author{Hao Cheng}
\address{School of Science, Jiangnan University, Jiangsu, Wuxi 214122, P. R. China}
\email{chenghao@jiangnan.edu.cn}

\author{Xiaoxiao Geng}

\thanks{}



\begin{abstract}
In this paper, we consider the prediction of the helium concentrations as function of a spatially variable source term perturbed by fractional Brownian motion. For the direct problem, we show that it is well-posed and has a unique mild solution under some conditions. For the inverse problem, the uniqueness and the instability are given. In the meanwhile, we determine the statistical properties of the source from the expectation and covariance of the final-time data $u(r,T)$. Finally, numerical implements are given to verify the effectiveness of the proposed reconstruction.
\end{abstract}

\maketitle

\section{Introduction}
The inverse source problem (ISP) of heat conduction is one of the important research subjects of inverse problems for partial differential equations. It has many scientific and industrial applications such as the detection and control of urban pollutants and marine pollution sources, the migration of groundwater, and the exploitation and control of unknown heat sources in oil wells. They have been extensively investigated in the literature \cite{3,4,6}. But the diffusion coefficient in the above research is constant, which is unreasonable in practical engineering applications. Meanwhile, we often encounter some special areas in practical problems. For instance, a steel storage tank yard for blast furnace steelmaking in the metallurgical field can be simplified as a cylindrical region; the helium diffusion of apatite in the laboratory is considered as a spherical area. In \cite{19}, wolf considered the model of helium production and diffusion in a spherical diffusion geometry:
\begin{align*}
 ^{4}He &= 8^{238}U(t)\left( e^{\lambda_{238}t}-1 \right) + 7^{235}U(t)\left( e^{\lambda_{235}t}-1 \right) + 6^{232}Th(t)\left( e^{\lambda_{232}t}-1 \right), \\
\frac{\partial^{4}He(r,t)}{\partial t} &= \frac{D(t)}{a^{2}}\left[ \frac{\partial^{2}{^{4}He(r,t)}}{\partial r^{2}} + \frac{2}{r}\frac{\partial^{4}He(r,t)}{\partial r}  \right] + 8 \lambda_{238}^{238}U(t) + 7 \lambda_{235}^{235}U(t) +6 \lambda_{232}^{232}Th(t),
\end{align*}
where $^{4}He(t)$, $U(t)$ and $Th(t)$ are amounts present at time t, $\lambda^{'}s$ are the decay constants, $a$ and $r$ are the radius and radial position of spherical diffusion domain, respectively, $D(t)$ is the time-dependent diffusion coefficient obeying an Arrhenius relationship such that:
$$D(t) = D_{0}exp \left[-\frac{E_{a}}{RT(t)} \right], $$
where $D_{0}$ is the diffusivity at infinite temperature and $E_{a}$ the activation energy measured in laboratory experiments, $R$ is the constant, and $T(t)$ is an arbitrary thermal history. Then the concentration of helium satisfies the heat conduction equation. In \cite{17}, Bao considered the model of helium production and diffusion in a spherical diffusion geometry and applied Tikhonov regularization and spectral truncation regularization methods to identify unknown source. Based on Bao's work, Zhang and Yan using a posteriori cutoff regularization and gave a posteriori convergence estimate in \cite{18}. Geng discussed the inverse source problem of the variable heat conduction equation on spherically symmetric domains by  using an iterative regularization method and obtained the H$\ddot{o}$lder type error estimates in \cite{7}.

However, the source is not being deterministic in the complex physical and engineering phenomena. The reasons for uncertainty can be roughly divided into the following three aspects \cite{9}: the first is the unpredictability of measurement environment, such as the change of temperature, wind speed and other objective conditions, the second is that the incomplete knowledge of the system will cause uncertainty in our experiment, the third is when the unit of measurement data is small, a very small error  can also produce great influence on the result.

Driven by these reasons, it is necessary to add random processes to the mathematical model. Since the solutions of stochastic differential equations are random functions, it is more important to study their statistical characteristics such as mean value, variance, and even higher order moments in many practical problems. In addition to instability, the ISP of stochastic equation is much more difficult than the deterministic problem due to randomness and uncertainty. Rensently, such problems have attracted more and more researchers' attention. In \cite{16}, Li analyzed the inverse scattering problem associated with the Schr$\ddot{o}$dinger system where both the potential and source terms are random and unknown. Bao reconstructed the statistical properties of random sources from boundary measurements of radiated random electric fields for the one-dimensional random Helmholtz equation in \cite{9}. Li showed that increasing stability can be obtained for the inverse problem by using only the Dirichlet boundary data with multi-frequencies for two or three dimensional Helmholtz equation in \cite{8}. In \cite{10}, Niu proved the existence and uniqueness of the direct problem and adopted the Tikhonov regularization method to solve the inverse stochastic source problem in diffusion equations.

The noise added in the above researches is the classical white noise $\dot{W}(t)$, which is the formal derivative of Brownian motion $W(t)$. Acyually in the more general case, the increments of the noise are not independent of each other, the source term is driven by a more general stochastic process. In \cite{20}, Li considered the random source problems for the time-harmonic acoustic and elastic wave equations in two and three dimensions, where the source term was assumed to be a microlocally isotropic generalized Gaussian random function. Li and Feng discussed the inverse random source problem for wave equation and time fractional diffusion equation, where the source was driven by a fractional Brownian motion in \cite{11, 12}.

In this paper, we focus on an ISP for the following variable coefficient heat equation in a spherically symmetric region and add a stochastic process in source term, which driven by the fractional Brownian motion
\begin{align}
\left\{ \begin{gathered}
u_{t}(r,t) = a(t) [u_{rr}(r,t)+\frac {2}{r}u(r,t)]+F(r,t), \qquad \hfill 0<r<R_{0}, 0<t<T, \\
u(r,0)=0,   \qquad \hfill 0 \leq r \leq R_{0} ,\\
u(R_{0},t)=0,  \qquad \hfill  0 \leq t \leq T ,\\
\lim_{r\rightarrow0} u (r,t)\quad bounded,  \qquad \hfill 0 \leq t \leq T ,
\end{gathered}
\right.
\end{align}
where $R_{0}$ is the radius of the sphere,  $a(t)$ is the diffusion coefficient that depends on the time variable, we suppose it satisfies
$$0 < a_{0} \leq a(t) \leq a_{1},$$
and the random source term $F(r,t)$ has the expression
 $$F(r,t)=f(r)h(t)+ g(r)\dot B^{H}(t).$$
Here $f(r)$ and $g(r)$ are deterministic functions with compact supports contained in spherical region, $h(t)$ is also a deterministic function, $ B ^{H}(t)$ is the fractional Brownian motion (fBm) with $H \in$ (0,1) which is called the Hurst index of fBm. The fBm is a widely used stochastic process that is particularly suited to model short- and long-range dependent phenomena, and anomalous diffusion in a variety of fields including physics, hydrology and financial mathematics, etc. In addition, the $\dot B^{H}(t)$ can be roughly understood as the derivative of $ B ^{H}(t)$ with respect to the time t.

The paper is organized as follows. In Section 2, we present some preliminary material. The direct problem satisfies some stability estimates are proved in section 3. In section 4, the uniqueness and instability is discussed for the inverse problem, which use the empirical expectation and correlation of the final-time data $u(r,T)$ to reconstruct $f$ and $\vert g \vert$ of the random source term $F$. Numerical examples to verify the effectiveness of our method are in provided Section 5. Finally, we give some concluding remarks in Section 6.

\section{Preliminaries}
\noindent Let $(\Omega, \mathcal{F}, \mathbb{P})$ be a complete probability space, where $\Omega$ is a sample space, $\mathcal{F}$ is a $\sigma$-algebra on $\Omega$, and $\mathbb{P}$ is probability measure on the measurable space $(\Omega, \mathcal{F})$. For a random variable $X$, we denote by $\mathbb{E}(X)$ and $\mathbb{V}(X)=\mathbb{E}(X - \mathbb{E}(X))^{2}=\mathbb{E}(X^{2})-(\mathbb{E}(X))^{2}$  the expectation and variance of $X$, respectively. For two random variables $X$ and $Y$ , $Cov(X,Y)=\mathbb{E}[(X - \mathbb{E}(X))(Y - \mathbb{E}(Y))]$ stands for the covariance of $X$ and $Y$. In the sequel, the dependence of random variables on the sample $\omega \in \Omega$ will be omitted unless it is necessary to avoid confusion.

The one-dimensional fBm $ B ^{H}(t)$ is a centered Gaussian process, which satisfies $B^{H}(0)= 0$ and is determined by the covariance function
$$R(t,s) = \mathbb{E}[B^{H}(t)B^{H}(s)]=\frac{1}{2}(t^{2H} + s^{2H}- \vert t-s \vert^{2H}),$$
for any $s,t \geq 0.$ In particular, if $H=\frac{1}{2}$, $B ^{H}(t)$ turns to be the standard Brownian motion $W(t)$, which has the covariance function $R(t,s)= t \land s.$

The increment of fBm satisfies
$$\mathbb{E}[(B ^{H}(t)-B^{H}(s))(B^{H}(s)-B^{H}(r))]=\frac{1}{2}[(t-r)^{2H}-(t-s)^{2H}-(s-r)^{2H}],$$
and
$$\mathbb{E}[(B^{H}(t)-B^{H}(s))^{2}]=(t-s)^{2H},$$
for any $0<r<s<t.$ It indicates that the increments of $B^{H}(t)$ in disjoint intervals are linearly dependent except for the case $H=\frac{1}{2}$, and the increments are stationary since their moments depend only on the length of the interval.

Based on the moment estimates and the Kolmogorov continuity criterion, it holds for any $\epsilon > 0$ and $s,t\in [0,T]$ that
$$\vert B^{H}(t)-B^{H}(s) \vert \leq C \vert t-s \vert^{H-\epsilon},$$
almost surely with constant $C$ depending on $\epsilon$ and $T.$ It is clear to note that $H$ represents the regularity of $B^{H}(t)$ and the trajectories of $B^{H}(t)$ are $(H-\epsilon )$- H$\ddot{o}$lder continuous.

The fBm $ B ^{H}(t)$ with $H \in (0,1)$ has a Wiener integral representation
$$ B ^{H}(t)= \int_{0}^{t}K_{H}(t,s)dW(s),$$
where $K_{H}$ is a square integrable kernel and $W$ is the standard Brownian motion.

For a fixed interval $[0,T]$, denote by $\mathcal{E}$ the space of step functions on $[0,T]$ and by $\mathcal{H}$ the closure of $\mathcal{E}$ with respect to the product
$$\langle \chi_{[0,t]}, \chi_{[0,s]} \rangle_{\mathcal{H}}= R(t,s),$$
where $\chi_{[0,t]}$ and $\chi_{[0,s]}$ are the characteristic functions. For $\psi(t), \phi(t) \in \mathcal{H}, $ we have that

(1) if $H \in (0, \frac{1}{2}),$
\begin{flalign}
&\ \quad \quad \mathbb{E}\bigg[ \int_{0}^{t}\psi(s)dB ^{H}(s) \int_{0}^{t}\phi(s)dB ^{H}(s)   \bigg] \nonumber \\ &
= \int_{0}^{t}\Big \{ c_{H}\big[ (\frac{t}{s})^{H-\frac{1}{2}}(t-s)^{H-\frac{1}{2}} -(H-\frac{1}{2})s^{\frac{1}{2}-H}\int_{s}^{t}u^{H-\frac{3}{2}}(u-s)^{h-\frac{1}{2}}du \big] \psi(s)\nonumber  \\  &
\quad + \int_{s}^{t}(\psi(u)-\psi(s))c_{H}(\frac{u}{s})^{H-\frac{1}{2}}(u-s)^{H-\frac{3}{2}}du \Big \} \nonumber \\ &
\quad \times \Big \{ c_{H}\big[ (\frac{t}{s})^{H-\frac{1}{2}}(t-s)^{H-\frac{1}{2}} -(H-\frac{1}{2})s^{\frac{1}{2}-H}\int_{s}^{t}u^{H-\frac{3}{2}}(u-s)^{h-\frac{1}{2}}du \big] \phi(s)\nonumber  \\  &
\quad + \int_{s}^{t}(\phi(u)-\phi(s))c_{H}(\frac{u}{s})^{H-\frac{1}{2}}(u-s)^{H-\frac{3}{2}}du \Big \} ds, &
\end{flalign}
where $c_{H}=(\frac{2H}{(1-2H)\beta (1-2H,H+\frac{1}{2})})^{\frac{1}{2}}$ and $\beta(p,q)= \int_{0}^{1}t^{p-1}(1-t)^{q-1}dt;$

(2) if $H=\frac{1}{2},$
\begin{align}
  \mathbb{E}\bigg[ \int_{0}^{t}\psi(s)dB ^{H}(s) \int_{0}^{t}\phi(s)dB ^{H}(s)   \bigg] &=\mathbb{E}\bigg[ \int_{0}^{t}\psi(s)dW(s) \int_{0}^{t}\phi(s)dW(s)   \bigg] \nonumber \\
&= \int_{0}^{t} \psi(s) \phi(s)ds;
\end{align}

(3) if $H \in (\frac{1}{2}, 1),$
\begin{align}
\mathbb{E}\bigg[ \int_{0}^{t}\psi(s)dB ^{H}(s) \int_{0}^{t}\phi(s)dB ^{H}(s)   \bigg] = \alpha_{H}\int_{0}^{t} \int_{0}^{t}\psi(s)\phi(s)\vert r-u \vert^{2H-2}dudr,
\end{align}
where $\alpha_{H}=H(2H-1).$

In this paper, $L^{2}([0,R_{0}];r^{2})$ represents the Hilbert space of the lebesgue measurable function $\phi(r)$ with weight $r^{2}$ on the interval $[0,R_{0}]$, $(\cdot,\cdot)$ and $\Vert \cdot \Vert $ represent the inner product and norm of the $L^{2}([0,R_{0}];r^{2})$ space respectively, which are defined as
$$(\phi(r),\omega_{n}(r)) = \int_{0}^{R_{0}} r^{2}\phi(r)\omega_{n}(r)dr,\quad \Vert \phi(r) \Vert = (\int_{0}^{R_{0}} r^{2}\vert \phi(r) \vert^{2})^{\frac{1}{2}}.$$

It is easy to note that the (1.1) is a classical direct problom if the diffusion confficient $a(t)$, the random source term $F(r,t)$, and the boundary conditions are known. If a solution to the direct problem exists, it must be unique. The only solution to the direct problem of the variable coefficient heat conduction equation without random terms has been discussed in \cite{7}. Similarly, we get the corresponding eigenvalues and eigenfunctions are
$$\lambda_{n}=\left( \frac{n\pi}{R_{0}} \right)^{2}, \quad \omega_{n}= \frac{\sqrt{2}n\pi}{\sqrt{ R_{0}^{3}}}\frac{\sin (n\pi r /R_{0})}{(n\pi r / R_{0})},\quad n=1,2,\cdots, $$
which in the interval $[0,R_{0}]$ is orthononal function system with weighted $r^{2},$ and completes in the class of functions squable over the interval $[0,R_{0}].$\\
{\bf Definition 1.} A stochastic process $u$ taking value in $L^{2}(D)$ is called a mild solution
\begin{align}
u(r,t)  &= \sum_{n=1}^{\infty}[f_{n}\int_{0}^{t}h(\tau)e^{-\lambda_{n}\int_{\tau}^{t}a(s)ds}d\tau+g_{n}\int_{0}^{t}e^{-\lambda_{n}\int_{\tau}^{t}a(s)ds}d B^{H}(\tau)]\omega_{n}(r) ,
\end{align}
which is equal to the following formula
\begin{align}
u_{n}(t)&=f_{n}\int_{0}^{t}h(\tau)e^{-\lambda_{n}\int_{\tau}^{t}a(s)ds}d\tau+g_{n}\int_{0}^{t}e^{-\lambda_{n}\int_{\tau}^{t}a(s)ds}d B^{H}(\tau) \nonumber \\
&:=  I_{n,1}(t) + I_{n,2}(t).
\end{align}
where
\begin{align}
f(r)=\sum_{n=0}^{\infty}f_{n}\omega_{n}(r), \quad f_{n}=\int_{0}^{R_{0}}r^{2} f(r)\omega_{n}(r)dr ; \nonumber \\
g(r)=\sum_{n=0}^{\infty}g_{n}\omega_{n}(r), \quad g_{n}=\int_{0}^{R_{0}}r^{2} g(r)\omega_{n}(r)dr. \nonumber
\end{align}

\section{The direct problem}
\noindent In this section, it is only necessary to address the stability since the existence and the uniqueness of the solution has already been considered \cite{13, 14, 15}. Given the following assumptions, we show the mild solution (2.8) of the initial-boundary value problem (1.1) is well-posedness .\\
{\bf Assumption 1.} Let $H \in (0,1)$ and $f, g\in L^{2}(D)$ with $\Vert g \Vert _{L^{2}(D)} \neq 0.$ Assume that $h\in L^{\infty}(0,T)$ has a positive lower bound $C(h)$ satisfying $h\geq C(h) >0$ $a.e.$ on the $(0,T).$

\begin{theorem}
Let Assumption 1 hold. Then the stochastic process u given in (2.8) satisfies
$$\mathbb{E}\bigg[ \Vert u \Vert_{L^{2}(D\times [0,T])}^{2}\bigg] \lesssim T^{3} \Vert f\Vert_{L^{2}(D)}^{2}\Vert h \Vert_{L^{\infty}(0,T)}^{2} +\frac{T^{2H+1}}{2H+1}\Vert g\Vert_{L^{2}(D)}^{2}.$$
\end{theorem}
 \noindent$Proof.$ It is easy to note that the mild solution (2.8) satisfies
\begin{align}
\Vert u(r,t) \Vert_{L^{2}(D)}^{2}&=\Big|\Big| \sum_{n=1}^{\infty}(I_{n,1}(t) + I_{n,2}(t))\omega_{n}(r)\Big|\Big|_{L^{2}(D)}^{2} \nonumber \\
&=\sum_{n=1}^{\infty}(I_{n,1}(t) + I_{n,2}(t))^{2} \nonumber \\
&\leq 2 \sum_{n=1}^{\infty}(I_{n,1}^{2}(t) + I_{n,2}^{2}(t)). \nonumber
\end{align}
Hence we have
\begin{align}
\mathbb{E}\bigg[ \Vert u \Vert_{L^{2}(D\times [0,T])}^{2}\bigg] &=\mathbb{E}\bigg[ \int_{0}^{T}\Vert u(r,t)\Vert_{L^{2}(D)}^{2} dt\bigg] \nonumber \\
&\lesssim \mathbb{E} \bigg[\int_{0}^{T} \bigg( \sum_{n=1}^{\infty}I_{n,1}^{2}(t) +\sum_{n=1}^{\infty} I_{n,2}^{2}(t)\bigg)dt \bigg]\nonumber \\
&=\int_{0}^{T} \bigg( \sum_{n=1}^{\infty}I_{n,1}^{2}(t)\bigg)dt +\mathbb{E}\bigg[\int_{0}^{T} \bigg( \sum_{n=1}^{\infty}I_{n,1}^{2}(t)\bigg)dt\bigg]\nonumber \\
&= \sum_{n=1}^{\infty}\Vert I_{n,1}(t) \Vert_{L^{2}(0,T)}^{2}+\int_{0}^{T}\bigg(\sum_{n=1}^{\infty}\mathbb{E}[I_{n,2}^{2}(t)]\bigg)dt\nonumber \\
&:=S_{1} +S_{2}.
\end{align}
We denote $a \lesssim b$ stands for $a \leq Cb,$ where $C>0$ is a constant. Obviously this notation can express the boundedness of $S_{1}$ and $S_{2}$, respectively.

For $S_{1}$, denoting $G_{n}(t)=e^{-a_{0}\lambda_{n}t}$, we obtain that
\begin{align}
\Vert I_{n,1}(t)\Vert_{L^{2}(0,T)}&=\Big|\Big| f_{n}\int_{0}^{t}h(\tau)e^{-\lambda_{n}\int_{\tau}^{t}a(s)ds}d\tau \Big|\Big|_{L^{2}(0,T)} \nonumber \\
&\leq \Big|\Big| f_{n}\int_{0}^{t}h(\tau)e^{-a_{0}\lambda_{n}(t-\tau)}d\tau \Big|\Big|_{L^{2}(0,T)} \nonumber \\
& = \left( \int_{0}^{T} \Big| f_{n}\int_{0}^{t}h(\tau)G_{n}(t-\tau)d\tau \Big| ^{2}dt\right)^{\frac{1}{2}} \nonumber \\
& \leq \left(\int_{0}^{T} \vert f_{n} \vert ^{2}\Vert h \Vert ^{2}_{L^{\infty}(0,T)}\Big( \int_{0}^{T}\vert G_{n}(\tau) \vert d\tau \Big)^{2} dt  \right)^{\frac{1}{2}} \nonumber \\
&\leq T^{\frac{1}{2}}\vert f_{n}\vert \Vert G_{n} \Vert_{L^{1}(0,T)}\Vert h \Vert_{L^{\infty}(0,T)}.
\end{align}
It is obvious that
\begin{align}
\Vert G_{n} \Vert_{L^{1}(0,T)} = \int_{0}^{T}e^{-a_{0}\lambda_{n}t}dt =\frac{1-e^{-a_{0}\lambda_{n}T}}{a_{0}\lambda_{n}} \leq T.
\end{align}
Combining (3.2) and (3.3) gives that
\begin{align}
S_{1}=\sum_{n=1}^{\infty}\Vert I_{n,1}(t) \Vert_{L^{2}(0,T)}^{2}\lesssim T \sum_{n=1}^{\infty}\vert f_{n}\vert^{2}\Vert G_{n} \Vert_{L^{1}(0,T)}^{2}\Vert h \Vert_{L^{\infty}(0,T)}^{2}=T^{3}\Vert h \Vert_{L^{\infty}(0,T)}^{2}\Vert f \Vert_{L^{2}(D)}^{2}.
\end{align}
For stochastic integral $S_{2}$, we note that
\begin{align*}
\mathbb{E}[I_{n,2}^{2}(t)]&=\mathbb{E} \bigg[ g_{n}^{2}\left( \int_{0}^{t}e^{-\lambda_{n}\int_{\tau}^{t}a(s)ds}d B^{H}(\tau)\right)^{2}\bigg]  \\
&=g_{n}^{2}\mathbb{E} \bigg[ \left( \int_{0}^{t}e^{-\lambda_{n}\int_{\tau}^{t}a(s)ds}d B^{H}(\tau)\right)^{2}\bigg] .
\end{align*}
Next, since the covariance operator of $B^{H}$ takes different forms in $H\in(0,\frac{1}{2})$, $H=\frac{1}{2}$, and $H\in(\frac{1}{2},1)$, we discuss the three cases separately.

For the case $H\in(0,\frac{1}{2})$, form (2.1) we have that
\begin{align}
&\mathbb{E}\bigg[ \left( \int_{0}^{t}e^{-\lambda_{n}\int_{\tau}^{t}a(s)ds}d B^{H}(\tau)\right)^{2}\bigg] \nonumber \\
&= \int_{0}^{t}\Big \{ c_{H}\left[ (\frac{t}{\tau})^{H-\frac{1}{2}}(t-\tau)^{H-\frac{1}{2}} -(H-\frac{1}{2})\tau^{\frac{1}{2}-H}\int_{\tau}^{t}u^{H-\frac{3}{2}}(u-\tau)^{H-\frac{1}{2}}du \right] \times e^{-\lambda_{n}\int_{\tau}^{t}a(s)ds}\nonumber  \\
&\quad + \int_{\tau}^{t}\left(e^{-\lambda_{n}\int_{u}^{t}a(s)ds}-e^{-\lambda_{n}\int_{\tau}^{t}a(s)ds}\right)c_{H}(\frac{u}{\tau})^{H-\frac{1}{2}}(u-\tau)^{H-\frac{3}{2}}du \Big \}^{2}d\tau  \nonumber  \\
&\lesssim \int_{0}^{t} \left[ (\frac{t}{\tau})^{H-\frac{1}{2}}(t-\tau)^{H-\frac{1}{2}}e^{-\lambda_{n}\int_{\tau}^{t}a(s)ds}\right]^{2}d\tau \nonumber \\
&\quad + \int_{0}^{t}\left[\tau^{\frac{1}{2}-H}\left( \int_{\tau}^{t}u^{H-\frac{3}{2}}(u-\tau)^{H-\frac{1}{2}}du \right)e^{-\lambda_{n}\int_{\tau}^{t}a(s)ds}\right]^{2}d\tau \nonumber \\
&\quad +\int_{0}^{t}\left[\int_{\tau}^{t}\bigg(e^{-\lambda_{n}\int_{u}^{t}a(s)ds}-e^{-\lambda_{n}\int_{\tau}^{t}a(s)ds}\bigg)(\frac{u}{\tau})^{H-\frac{1}{2}}(u-\tau)^{H-\frac{3}{2}}du\right]^{2}d\tau \nonumber \\
&:= I_{1}(t) + I_{2}(t) + I_{3}(t).
\end{align}
A direct calculation yields
\begin{align}
I_{1}(t) &=\int_{0}^{t} (\frac{t}{\tau})^{2H-1}(t-\tau)^{2H-1} e^{-2\lambda_{n}\int_{\tau}^{t}a(s)ds} d\tau \nonumber \\
&\leq \int_{0}^{t} (\frac{t}{\tau})^{2H-1}(t-\tau)^{2H-1} d\tau \nonumber \\
&\leq \int_{0}^{t} (t-\tau)^{2H-1} d\tau = \frac{t^{2H}}{2H}.
\end{align}
Similarly, we have
\begin{align}
I_{2}(t)&=\int_{0}^{t}\tau^{1-2H}\left( \int_{\tau}^{t}u^{H-\frac{3}{2}}(u-\tau)^{H-\frac{1}{2}}du \right)^{2}e^{-2\lambda_{n}\int_{\tau}^{t}a(s)ds} d\tau \nonumber \\
&\leq \int_{0}^{t}\tau^{1-2H}\left( \int_{\tau}^{t}u^{H-\frac{3}{2}}(u-\tau)^{H-\frac{1}{2}}du \right)^{2}d\tau.
\end{align}
From the binomial expansion and the Raabe discriminant, it is easy to know the following series is absolutely convergent. Therefore, we can deduce that
\begin{align}
\int_{\tau}^{t}u^{H-\frac{3}{2}}(u-\tau)^{H-\frac{1}{2}}du
&=\int_{\tau}^{t} u^{2H-2}\left[ \sum_{n=0}^{\infty}{H-\frac{1}{2} \choose n}(-\frac{\tau}{u})^{n} \right] du \nonumber \\
&=\sum_{n=0}^{\infty}{H-\frac{1}{2} \choose n}(-1)^{n}\tau^{n}\left( \frac{\tau^{2H-1-n}-t^{2H-1-n}}{n+1-2H} \right)\nonumber \\
&\leq (\tau^{2H-1}-t^{2H-1})\sum_{n=0}^{\infty}{H-\frac{1}{2} \choose n} \frac{(-1)^{n}}{n+1-2H}\nonumber \\
&\lesssim \tau^{2H-1}-t^{2H-1}.
\end{align}
Combining (3.6) and (3.7), we obtain that
\begin{align}
I_{2}(t)& \lesssim \int_{0}^{t}\tau^{1-2H}(\tau^{2H-1}-t^{2H-1})^{2} d\tau \nonumber \\
&\lesssim \int_{0}^{t}\tau^{1-2H}(t^{4H-2}+\tau^{4H-2}) d\tau \nonumber \\
&\lesssim \int_{0}^{t}\tau^{2H-1} d\tau =\frac{t^{2H}}{2H}.
\end{align}
From simple inequalities $e^{-r} \leq r^{-\frac{1}{2}}$ and $1-e^{-r}\leq r^{\frac{1}{2}}(r>0)$, We have
\begin{align}
I_{3}(t) &= \int_{0}^{t}\left[\int_{\tau}^{t}\bigg(e^{-\lambda_{n}\int_{u}^{t}a(s)ds}-e^{-\lambda_{n}\int_{\tau}^{t}a(s)ds}\bigg)(\frac{u}{\tau})^{H-\frac{1}{2}}(u-\tau)^{H-\frac{3}{2}}du\right]^{2}d\tau \nonumber \\
&= \int_{0}^{t}\left[\int_{\tau}^{t}e^{-\lambda_{n}\int_{u}^{t}a(s)ds}\bigg(1-e^{-\lambda_{n}\int_{\tau}^{u}a(s)ds}\bigg)(\frac{u}{\tau})^{H-\frac{1}{2}}(u-\tau)^{H-\frac{3}{2}}du\right]^{2}d\tau \nonumber \\
&= \int_{0}^{t}\left[\int_{\tau}^{t}e^{-\lambda_{n}a_{0}(t-u)}\bigg(1-e^{-\lambda_{n}a_{1}(u-\tau)}\bigg)(\frac{u}{\tau})^{H-\frac{1}{2}}(u-\tau)^{H-\frac{3}{2}}du\right]^{2}d\tau \nonumber \\
&\lesssim \int_{0}^{t}\left[\int_{\tau}^{t}(t-u)^{-\frac{1}{2}}(u-\tau)^{\frac{1}{2}}(\frac{u}{\tau})^{H-\frac{1}{2}}(u-\tau)^{H-\frac{3}{2}}du\right]^{2}d\tau \nonumber \\
& \lesssim \int_{0}^{t}\left[\int_{\tau}^{t}(t-u)^{-\frac{1}{2}}(\frac{u}{\tau})^{H-\frac{1}{2}}(u-\tau)^{H-1}du\right]^{2}d\tau, \nonumber
\end{align}
which is integratable under the condition $H\in(0,\frac{1}{2})$. Since $0<\tau<u<t$, we have $(\frac{u}{\tau})^{H-\frac{1}{2}} <1$. From Raabe discriminant, we make $q=u-\tau$ to replace variable and obtain that
\begin{align}
I_{3}(t) &\lesssim \int_{0}^{t}\left[\int_{0}^{t-\tau}(t-\tau-q)^{-\frac{1}{2}}(q)^{H-1}dq\right]^{2}d\tau \nonumber \\
&=\int_{0}^{t}\left[\int_{0}^{t-\tau} \sum_{n=0}^{\infty}\binom{-\frac{1}{2}}{n}(-1)^{n}\left(\frac{t-\tau}{q}\right)^{-\frac{1}{2}-n} (q)^{H-\frac{3}{2}}dq\right]^{2}d\tau \nonumber \\
&=\int_{0}^{t}\left[\int_{0}^{t-\tau} \sum_{n=0}^{\infty}\binom{-\frac{1}{2}}{n}(-1)^{n}\left(t-\tau \right)^{-\frac{1}{2}-n}  (q)^{n+H-1}dq\right]^{2}d\tau \nonumber \\
&= \left[\sum_{n=0}^{\infty}\binom{-\frac{1}{2}}{n}\frac{(-1)^{n}}{n+H}       \right]^{2} \int_{0}^{t}(t-\tau)^{2H-1}d\tau \nonumber \\
& \lesssim t^{2H}.
\end{align}

Substituting (3.6), (3.9) and (3.10) into (3.5), we can get when $H\in(0,\frac{1}{2}),$ that
$$\mathbb{E}\bigg[ \left( \int_{0}^{t}e^{-\lambda_{n}\int_{\tau}^{t}a(s)ds}d B^{H}(\tau)\right)^{2}\bigg] \lesssim t^{2H}.$$

For the case $H=\frac{1}{2},$ it follows from It$\hat{o}$ isometry (2.2) that
\begin{align*}
\mathbb{E}\bigg[ \left( \int_{0}^{t}e^{-\lambda_{n}\int_{\tau}^{t}a(s)ds}d B^{\frac{1}{2}}(\tau)\right)^{2}\bigg] &= \int_{0}^{t}e^{-2\lambda_{n}\int_{\tau}^{t}a(s)ds}d \tau \\
& \leq \int_{0}^{t}e^{-2\lambda_{n}a_{0}(t-\tau)}d \tau \\
& = \frac{1-e^{-2\lambda_{n}a_{0}t}}{2\lambda_{n}a_{0}} \leq t.
\end{align*}

For the case $H\in(\frac{1}{2},1)$, we have from (2.3) that
\begin{align}
\mathbb{E}\bigg[ \left( \int_{0}^{t}e^{-\lambda_{n}\int_{\tau}^{t}a(s)ds}d B^{H}(\tau)\right)^{2}\bigg]&=\alpha_{H}\int_{0}^{t}\int_{0}^{t}e^{-\lambda_{n}\int_{r}^{t}a(s)ds}e^{-\lambda_{n}\int_{u}^{t}a(s)ds}\vert r-u \vert^{2H-2} dudr \nonumber \\
&\leq \alpha_{H}\int_{0}^{t}\int_{0}^{t}\vert r-u \vert^{2H-2}dudr \nonumber \\
&\lesssim \int_{0}^{t} \left( \frac{(t-r)^{2H-1}}{2H-1}+\frac{r^{2H-1}}{2H-1}\right)dr \nonumber \\
&\lesssim t^{2H}. \nonumber
\end{align}

Therefore, for any $H\in(0,1),$ it holds that
\begin{align*}
\mathbb{E}\bigg[ \left( \int_{0}^{t}e^{-\lambda_{n}\int_{\tau}^{t}a(s)ds}d B^{H}(\tau)\right)^{2}\bigg]\lesssim t^{2H}.
\end{align*}
Hence
\begin{align}
S_{2}=\int_{0}^{T}\bigg(\sum_{n=1}^{\infty}\mathbb{E}[I_{n,2}^{2}(t)]\bigg)dt &\lesssim \int_{0}^{T}\left( \sum_{n=1}^{\infty} g_{n}^{2}\mathbb{E}\bigg[ \left( \int_{0}^{t}e^{-\lambda_{n}\int_{\tau}^{t}a(s)ds}d B^{H}(\tau)\right)^{2}\bigg] \right) dt \nonumber \\
&\lesssim \int_{0}^{T} \sum_{n=1}^{\infty} g_{n}^{2} t^{2H}dt = \frac{T^{2H+1}}{2H+1}\Vert g\Vert_{L^{2}(D)}^{2}.
\end{align}
From (3.4) and (3.11) we obtain the following conclution
$$\mathbb{E}\bigg[ \Vert u \Vert_{L^{2}(D\times [0,T])}^{2}\bigg] \lesssim T^{3} \Vert f\Vert_{L^{2}(D)}^{2}\Vert h \Vert_{L^{\infty}(0,T)}^{2} +\frac{T^{2H+1}}{2H+1}\Vert g\Vert_{L^{2}(D)}^{2},$$
which implies the stability estimate for the mild solution (2.8).

\section{The inverse problem}
\noindent In this section, we attempt to reconstruct $f$ and $\left| g\right|$ from the empirical expectation and correlation of the final time data $u(r,T)$. From (2.1)-(2.3) and(2.5), it follows that the final time expectation and variance can be formulated as
\begin{align}
\mathbb{E}(u_{n}(T))&=f_{n}\int_{0}^{T}h(\tau)e^{-\lambda_{n}\int_{\tau}^{T}a(s)ds}d\tau, \\
\mathbb{V}(u_{n}(T))&=g_{n}^{2}\mathbb{E}\left[ \left( \int_{0}^{T}e^{-\lambda_{n}\int_{\tau}^{T}a(s)ds}d B^{H}(\tau) \right) ^{2}\right],
\end{align}
and
\begin{align}
Cov(u_{m}(T),u_{n}(T))=g_{m}g_{n}\mathbb{E}\left[ \int_{0}^{T}e^{-\lambda_{m}\int_{\tau}^{T}a(s)ds}d B^{H}(\tau)\int_{0}^{T}e^{-\lambda_{n}\int_{\tau}^{T}a(s)ds}d B^{H}(\tau) \right].
\end{align}
Below, we consider the uniqueness and the instability of the inverse problem, separately.\\
4.1 {\bf Uniqueness.}
\begin{theorem}
Let assumption 1 hold. Then the source terms $f$ and $g$ up to sign, i.e. $\left| g\right|$  can be uniquely determined by the data set
$$\{ \mathbb{E}(u_{n}(T)),Cov(u_{m}(T),u_{n}(T)): m,n \in \mathbb{N} \}.$$
\end{theorem}
\noindent $Proof.$ Firstly, we need to prove for each fixed $n \in \mathbb{N}$ there exists a constant $C_{1} > 0$ such that
\begin{align}
\int_{0}^{T}h(\tau)e^{-\lambda_{n}\int_{\tau}^{T}a(s)ds}d\tau \geq C_{1} >0.
\end{align}
setting $\tilde{\tau} = T - \tau,$ from assumption 1 we have that
\begin{align*}
\int_{0}^{T}h(\tau)e^{-\lambda_{n}\int_{\tau}^{T}a(s)ds}d\tau
&= \int_{0}^{T}h(T-\tilde{\tau})e^{-\lambda_{n}\int_{T-\tilde{\tau}}^{T}a(s)ds}d\tilde{\tau} \\
& \geq \int_{0}^{T}h(T-\tilde{\tau})e^{-\lambda_{n}a_{1}\tilde{\tau}} d\tilde{\tau} \\
& \geq C(h) \int_{0}^{T} e^{-\lambda_{n}a_{1}\tilde{\tau}}d\tilde{\tau}\\
&= C(h)\frac{1-e^{-\lambda_{n}a_{1}T}}{\lambda_{n}a_{1}}:= C_{1}>0.
\end{align*}

Secondly, we certificate for each fixed $n, m \in \mathbb{N}$ there exists another constant $C_{2} > 0$ such that
\begin{align}
\mathbb{E}\left[ \int_{0}^{T}e^{-\lambda_{m}\int_{\tau}^{T}a(s)ds}d B^{H}(\tau)\int_{0}^{T}e^{-\lambda_{n}\int_{\tau}^{T}a(s)ds}d B^{H}(\tau) \right] \geq C_{2} >0.
\end{align}
We denote $\phi_{n}(\tau)=e^{-\lambda_{n}\int_{\tau}^{T}a(s)ds}$ and
\begin{align*}
&\quad \mathbb{E}\left[ \int_{0}^{T}e^{-\lambda_{m}\int_{\tau}^{T}a(s)ds}d B^{H}(\tau)\int_{0}^{T}e^{-\lambda_{n}\int_{\tau}^{T}a(s)ds}d B^{H}(\tau) \right] \\
&= \mathbb{E}\left[ \int_{0}^{T}\phi_{m}(\tau)d B^{H}(\tau)\int_{0}^{T}\phi_{n}(\tau)d B^{H}(\tau) \right] :=E_{mn}.
\end{align*}
Next, we estimate $E_{mn}$ for $H\in(0,\frac{1}{2})$, $H=\frac{1}{2}$, and $H\in(\frac{1}{2},1)$, separately.

For $H\in(0,\frac{1}{2}),$ we first note that  the square integrable kernel $K_{H}$ has the following form
$$K_{H}(T,s)=c_{H}\left[ (\frac{T}{s})^{H-\frac{1}{2}}(T-s)^{H-\frac{1}{2}} -(H-\frac{1}{2})s^{\frac{1}{2}-H}\int_{s}^{t}u^{H-\frac{3}{2}}(u-s)^{h-\frac{1}{2}}du \right], $$
Since $H\in (0,\frac{1}{2}),$ $0<\tau<T$, it is obvious that $K_{H}(T,\tau) >0 $. From the mean value theorem, it follows that
\begin{align*}
&\quad \int_{s}^{T}(\phi_{n}(u)-\phi_{n}(\tau))c_{H}(\frac{u}{\tau})^{H-\frac{1}{2}}(u-\tau)^{H-\frac{3}{2}}du\\
&= c_{H}\int_{\tau}^{T} \phi_{n}^{'}(u_{k}^{*})(\frac{u}{\tau})^{H-\frac{1}{2}}(u-\tau)^{H-\frac{1}{2}}du \qquad \quad (\tau<u_{k}^{*}<u<T) \\
&= M_{H}(\tau)\phi_{n}^{'}(u_{k}^{*}),
\end{align*}
where
$$M_{H}(\tau)=c_{H}\int_{\tau}^{T}(\frac{u}{\tau})^{H-\frac{1}{2}}(u-\tau)^{H-\frac{1}{2}}du > 0.$$
Denote by $\langle \cdot , \cdot \rangle$ the inner product on $L^{2}(0,T)$, which given by
$$\langle \psi(\tau),\phi(\tau)\rangle =\int_{0}^{T}\psi(\tau)\phi(\tau)d\tau .$$
We have from (2.1) that
\begin{align*}
E_{mn} &= \langle K_{H}(T,\tau)\phi_{m}(\tau)+M_{H}(\tau)\phi_{m}^{'}(u_{k}^{*}),K_{H}(T,\tau)\phi_{n}(\tau)+M_{H}(\tau)\phi_{n}^{'}(u_{k}^{*}) \rangle_{L^{2}(0,T)}\\
&= \langle K_{H}(T,\tau)\phi_{m}(\tau),K_{H}(T,\tau)\phi_{n}(\tau)\rangle_{L^{2}(0,T)} + \langle M_{H}(\tau)\phi_{m}^{'}(u_{k}^{*}),M_{H}(\tau)\phi_{n}^{'}(u_{k}^{*}) \rangle_{L^{2}(0,T)} \\
 & \quad +\langle K_{H}(T,\tau)\phi_{m}(\tau),M_{H}(\tau)\phi_{n}^{'}(u_{k}^{*}) \rangle_{L^{2}(0,T)} +\langle M_{H}(\tau)\phi_{m}^{'}(u_{k}^{*}),K_{H}(T,\tau)\phi_{n}(\tau)\rangle_{L^{2}(0,T)}.
\end{align*}
We consider the first term
\begin{align*}
&\quad \langle K_{H}(T,\tau)\phi_{m}(\tau),K_{H}(T,\tau)\phi_{n}(\tau)\rangle_{L^{2}(0,T)} \\
&= \int_{0}^{T}K_{H}^{2}(T,\tau)\phi_{m}(\tau)\phi_{n}(\tau) d\tau \\
& \geq e^{-(\lambda_{m}+\lambda_{n})a_{1}T}\int_{0}^{T}K_{H}^{2}(T,\tau)d\tau :=\tilde{c}_{1} >0.
\end{align*}
We notice that $\phi_{n}(\tau) >0$ and $\phi_{n}(\tau)$ is a monotonically increasing function as $\phi_{n}^{'}(\tau)=\lambda_{n}a(\tau)\phi_{n}(\tau)>0$, $\phi_{n}^{'}(\tau)$ is a monotonically increasing function as $\phi_{n}^{''}(\tau)=\lambda_{n}\phi_{n}(\tau)+\lambda_{n}^{2}a^{2}(\tau)\phi_{n}(\tau) >0$. Accordingly we have $\phi_{n}(\tau) \geq\phi_{n}(0) >0$ and $\phi_{n}^{'}(\tau)\geq \phi_{n}^{'}(0) > 0$, for $0<\tau<T$. Hence
$$\langle M_{H}(\tau)\phi_{m}^{'}(u_{k}^{*}),M_{H}(\tau)\phi_{n}^{'}(u_{k}^{*}) \rangle_{L^{2}(0,T)} \geq \phi_{m}^{'}(0)\phi_{n}^{'}(0)\int_{0}^{T}M_{H}^{2}(\tau)d\tau:=\tilde{c}_{2} >0. $$
Similarly, we have
$$\langle K_{H}(T,\tau)\phi_{m}(\tau),M_{H}(\tau)\phi_{n}^{'}(u_{k}^{*}) \rangle_{L^{2}(0,T)}\geq \phi_{m}(0)\phi_{n}^{'}(0) \int_{0}^{T}K_{H}(T,\tau)M_{H}(\tau)d\tau :=\tilde{c}_{3} >0, $$
and
$$\langle M_{H}(\tau)\phi_{m}^{'}(u_{k}^{*}),K_{H}(T,\tau)\phi_{n}(\tau)\rangle_{L^{2}(0,T)}\geq \phi_{m}^{'}(0)\phi_{n}(0) \int_{0}^{T}K_{H}(T,\tau)M_{H}(\tau)d\tau :=\tilde{c}_{4} >0.$$
From the above estimates we conclude that
$$E_{mn} \geq \sum_{j=1}^{4}\tilde{c}_{j}:= C_{2} >0.$$

For $H= \frac{1}{2},$ from (2.2), we have
\begin{align*}
E_{mn} = \mathbb{E}\left[ \int_{0}^{T}e^{-(\lambda_{m}+\lambda_{n})\int_{\tau}^{T}a(s)ds}d \tau \right]
 \geq \int_{0}^{T}e^{-(\lambda_{m}+\lambda_{n})a_{1}T}d\tau > 0.
\end{align*}

For $H \in (\frac{1}{2},1),$ from (2.3), we have
\begin{align*}
E_{mn} &= \alpha_{H}\int_{0}^{T}\int_{0}^{T}e^{-\lambda_{m}\int_{r}^{T}a(s)ds}e^{-\lambda_{m}\int_{u}^{T}a(s)ds}\vert r-u \vert^{2H-2}dudr\\
& \geq \alpha_{H}\int_{0}^{T}\int_{0}^{T}e^{-\lambda_{m}a_{1}T}e^{-\lambda_{n}a_{1}T}\vert r-u \vert^{2H-2}dudr\\
& = \alpha_{H}e^{-(\lambda_{m}+\lambda_{n})a_{1}T}\int_{0}^{T}\int_{0}^{T}\vert r-u \vert^{2H-2}dudr \\
& =T^{2H}e^{-(\lambda_{m}+\lambda_{n})a_{1}T} >0.
\end{align*}

Finally, we have from $f,g \in L^{2}(D)$ that
$$f(r)=\sum_{n=1}^{\infty}f_{k}\omega_{n}(r),\qquad g(r)=\sum_{n=1}^{\infty}g_{k}\omega_{n}(r) ,$$
which gives that
$$g^{2}(x)=\left( \sum_{m=1}^{\infty}g_{n}\omega_{m}(r) \right)\left( \sum_{n=1}^{\infty}g_{n}\omega_{n}(r) \right)= \sum_{m,n\in \mathbb{N}}g_{n}g_{n}\omega_{m}(r)\omega_{n}(r).$$

Combining (4.1)-(4.5), we can obtain the uniqueness of the inverse problem, which completes the proof.\\
4.1 {\bf Instability.}

\begin{theorem}
The inverse problem is ill-posedness to reconstruct the source terms f and $\pm g.$ Moreover, if $h \in L^{\infty}(0,T)$ is a nonnegative function and its support has a positive measure, then the following estimates hold
$$\Big \vert \int_{0}^{T}h(\tau)e^{-\lambda_{n}\int_{\tau}^{T}a(s)ds}d\tau \Big \vert \lesssim \frac{1}{\lambda_{n}},$$
$$\mathbb{E}\left[ \left( \int_{0}^{T}e^{-\lambda_{n}\int_{\tau}^{T}a(s)ds}d B^{H}(\tau) \right) ^{2}\right] \lesssim \frac{1}{\lambda_{n}^{p}},$$
where p=min\{2H,1\}.
\end{theorem}

\noindent$Proof.$ First, it is clear to note that
$$\int_{0}^{T}e^{-\lambda_{n}\int_{\tau}^{T}a(s)ds}d\tau \leq  \int_{0}^{T}e^{-\lambda_{n}a_{0}(T-\tau)}d\tau =\frac{1-e^{-\lambda_{n}a_{0}T}}{\lambda_{n}a_{0}} <\frac{1}{\lambda_{n}a_{0}},$$
\begin{align}
\Big \vert \int_{0}^{T}h(\tau)e^{-\lambda_{n}\int_{\tau}^{T}a(s)ds}d\tau \Big \vert  \lesssim \frac{1}{\lambda_{n}}  \quad \rightarrow 0\quad as \quad n\rightarrow \infty,
\end{align}
which shows that it is unstable to recover $f$ due to (4.1).

Second,we discuss the instability of recovering $g_{n}^{2}$, which is equivalent to the instability of recovering $\vert g \vert .$ Similarly, we present the three different cases $H\in (0,\frac{1}{2})$, $H=\frac{1}{2}$ and $H \in (\frac{1}{2},1)$, separately.

For the case $H\in(0,\frac{1}{2})$, we consider (3.6) with $t=T$ and the estimate of $I_{j}(T),$ $j=1,2,3.$ For $I_{1}(T),$ a simple calculation yields
\begin{align*}
I_{1}(T) &= \int_{0}^{T} (\frac{\tau}{T})^{1-2H}(T-\tau)^{2H-1}e^{-2\lambda_{n}\int_{\tau}^{T}a(s)ds} d\tau \\
&\leq e^{-2\lambda_{n}a_{0}T} \int_{0}^{T} (\frac{\tau}{T})^{1-2H}(T-\tau)^{2H-1} d\tau \\
&\leq e^{-2\lambda_{n}a_{0}T} \int_{0}^{T} (T-\tau)^{2H-1} d\tau \\
&\lesssim \frac{1}{\lambda_{n}}\frac{T^{2H-1}}{2H}.
\end{align*}
For $I_{2}(T),$ we have
\begin{align*}
I_{2}(T)&= \int_{0}^{T}\tau^{1-2H}\left( \int_{\tau}^{T}u^{H-\frac{3}{2}}(u-\tau)^{H-\frac{1}{2}}du \right)^{2}e^{-2\lambda_{n}\int_{\tau}^{T}a(s)ds}d\tau \\
& \leq e^{-2\lambda_{n}a_{0}T} \int_{0}^{T}\tau^{1-2H}\left( \int_{\tau}^{T}u^{H-\frac{3}{2}}(u-\tau)^{H-\frac{1}{2}}du \right)^{2} d\tau \\
&\lesssim e^{-2\lambda_{n}a_{0}T} \int_{0}^{T}\tau^{1-2H}(T^{2H-1}-\tau^{2H-1})^{2} d\tau \\
& \lesssim e^{-2\lambda_{n}a_{0}T}\int_{0}^{T}\tau^{1-2H}(T^{4H-2}+\tau^{4H-2}) d\tau \\
& \lesssim e^{-2\lambda_{n}a_{0}T}\int_{0}^{T}\tau^{2H-1}d\tau \lesssim \frac{1}{\lambda_{n}}\frac{T^{2H-1}}{2H}.
\end{align*}
For $I_{3}(T)$, we divide it as follows
\begin{align*}
I_{3}(T) &= \int_{0}^{T}\left[\int_{\tau}^{T}\bigg(e^{-\lambda_{n}\int_{u}^{T}a(s)ds}-e^{-\lambda_{n}\int_{\tau}^{T}a(s)ds}\bigg)(\frac{u}{\tau})^{H-\frac{1}{2}}(u-\tau)^{H-\frac{3}{2}}du\right]^{2}d\tau \\
&= \int_{0}^{t_{*}}\left[ \int_{\tau}^{t_{*}}\bigg(e^{-\lambda_{n}\int_{u}^{T}a(s)ds}-e^{-\lambda_{n}\int_{\tau}^{T}a(s)ds}\bigg)(\frac{u}{\tau})^{H-\frac{1}{2}}(u-\tau)^{H-\frac{3}{2}}du  \right]^{2}d\tau \\
&\quad +\int_{0}^{t_{*}}\left[ \int_{t_{*}}^{T}\bigg(e^{-\lambda_{n}\int_{u}^{T}a(s)ds}-e^{-\lambda_{n}\int_{\tau}^{T}a(s)ds}\bigg)(\frac{u}{\tau})^{H-\frac{1}{2}}(u-\tau)^{H-\frac{3}{2}}du  \right]^{2}d\tau \\
&\quad +\int_{t_{*}}^{T}\left[ \int_{\tau}^{\tau+t_{*}}\bigg(e^{-\lambda_{n}\int_{u}^{T}a(s)ds}-e^{-\lambda_{n}\int_{\tau}^{T}a(s)ds}\bigg)(\frac{u}{\tau})^{H-\frac{1}{2}}(u-\tau)^{H-\frac{3}{2}}du  \right]^{2}d\tau \\
&\quad +\int_{t_{*}}^{T}\left[ \int_{\tau+t_{*}}^{T}\bigg(e^{-\lambda_{n}\int_{u}^{T}a(s)ds}-e^{-\lambda_{n}\int_{\tau}^{T}a(s)ds}\bigg)(\frac{u}{\tau})^{H-\frac{1}{2}}(u-\tau)^{H-\frac{3}{2}}du  \right]^{2}d\tau \\
&:= K_{1}+K_{2}+K_{3}+K_{4},
\end{align*}
where $t_{*}\in(0,T)$ will be chosen later. About $K_{1}$ and $K_{2},$ from the differential mean value theorem we obtain that
\begin{align*}
K_{1} & \leq \lambda_{n}a(u_{*})e^{-\lambda_{n}\int_{u_{*}}^{T}a(s)ds}\int_{0}^{t_{*}}\left[ \int_{\tau}^{t_{*}}(\frac{u}{\tau})^{H-\frac{1}{2}}(u-\tau)^{H-\frac{1}{2}}du  \right]^{2}d\tau \\
& \leq \frac{a(u_{*})}{\int_{u_{*}}^{T}a(s)ds}\int_{0}^{t_{*}}\left[ \int_{\tau}^{t_{*}}(u-\tau)^{H-\frac{1}{2}}du  \right]^{2}d\tau \\
& \lesssim \int_{0}^{t_{*}}\left[ \int_{\tau}^{t_{*}}(u-\tau)^{H-\frac{1}{2}}du  \right]^{2}d\tau \\
&\lesssim \int_{0}^{t_{*}}(t_{*}-\tau)^{2H-1}d\tau \lesssim t_{*}^{2H+2},
\end{align*}
similarly, it is easy to get
\begin{align*}
K_{3} & \lesssim \int_{t_{*}}^{T}\left[ \int_{\tau}^{\tau+t_{*}}(u-\tau)^{H-\frac{1}{2}}du  \right]^{2}d\tau \\
&\lesssim \int_{t_{*}}^{T}t_{*}^{2H+1} d\tau \lesssim t_{*}^{2H+1}.
\end{align*}
About $K_{2}$ and $K_{4},$ we obtain from straightforward calculations that
\begin{align*}
K_{2} & \leq \int_{0}^{t_{*}}\left[ \int_{t_{*}}^{T}e^{-\lambda_{n}a_{0}(T-u)}(\frac{u}{\tau})^{H-\frac{1}{2}}(u-\tau)^{H-\frac{3}{2}}du  \right]^{2}d\tau \\
& \lesssim \frac{1}{\lambda_{n}} \int_{0}^{t_{*}}\left[ \int_{t_{*}}^{T}(u-\tau)^{H-\frac{3}{2}}du  \right]^{2}d\tau \\
& \lesssim \frac{1}{\lambda_{n}} \int_{0}^{t_{*}}\left[ (T-\tau)^{2H-1}-(t_{*}-\tau)^{2H-1} \right]d\tau \\
& \lesssim \frac{1}{\lambda_{n}} \left[T^{2H}-(T-t_{*})^{2H}+ t_{*}^{2H} \right]\\
&\lesssim \frac{1}{\lambda_{n}} t_{*}^{2H},
\end{align*}
and
\begin{align*}
K_{4} & \leq \int_{t_{*}}^{T}\left[ \int_{\tau+t_{*}}^{T}e^{-\lambda_{n}a_{0}(T-u)}(\frac{u}{\tau})^{H-\frac{1}{2}}(u-\tau)^{H-\frac{3}{2}}du  \right]^{2}d\tau \\
&\lesssim \frac{1}{\lambda_{n}}\int_{t_{*}}^{T}\left[ \int_{\tau+t_{*}}^{T}(u-\tau)^{H-\frac{3}{2}}du  \right]^{2}d\tau \\
&\lesssim \frac{1}{\lambda_{n}}\int_{t_{*}}^{T}\left[(T-\tau)^{2H-1}+ t_{*}^{2H-1} \right]d\tau \\
&\lesssim \frac{1}{\lambda_{n}}(T-t_{*})^{2H}+ \frac{1}{\lambda_{n}}t_{*}^{2H-1}  \\
&\lesssim \frac{1}{\lambda_{n}}t_{*}^{2H} + \frac{1}{\lambda_{n}}t_{*}^{2H-1}.
\end{align*}
Combining the above estimates about $K_{1}-K_{4}$, and choosing $t_{*}=1/\lambda_{n},$ we get
$$I_{3}(T)\lesssim t_{*}^{2H+2}+\frac{1}{\lambda_{n}} t_{*}^{2H} +t_{*}^{2H+1} +\frac{1}{\lambda_{n}}t_{*}^{2H} + \frac{1}{\lambda_{n}}t_{*}^{2H-1} \lesssim \frac{1}{\lambda_{n}^{2H}}.$$
From the above estimates for $I_{1}, I_{2}$ and $I_{3}$, we get for $H\in(0,\frac{1}{2})$ that
\begin{align*}
\mathbb{E}\left[ \left( \int_{0}^{T}e^{-\lambda_{n}\int_{\tau}^{T}a(s)ds}d B^{H}(\tau) \right) ^{2}\right] \lesssim \frac{1}{\lambda_{n}}\frac{T^{2H-1}}{H} + \frac{1}{\lambda_{n}^{2H}}.
\end{align*}

For the case $H=\frac{1}{2}$, we have
\begin{align*}\mathbb{E}\left[ \left( \int_{0}^{T}e^{-\lambda_{n}\int_{\tau}^{T}a(s)ds}d B^{\frac{1}{2}}(\tau) \right) ^{2}\right] &= \int_{0}^{T}e^{-2\lambda_{n}\int_{\tau}^{T}a(s)ds}d\tau \leq \int_{0}^{T} e^{-2\lambda_{n}a_{0}(T-\tau)}d\tau \\
& = \frac{1-e^{-2\lambda_{n}a_{0}T}}{2\lambda_{n}a_{0}} \lesssim \frac{1}{\lambda_{n}}.
\end{align*}

For the case $H\in(\frac{1}{2},1)$, we have
\begin{align*}
& \mathbb{E}\left[ \left( \int_{0}^{T}e^{-\lambda_{n}\int_{\tau}^{T}a(s)ds}d B^{H}(\tau) \right) ^{2}\right] \\
&= \alpha_{H}\int_{0}^{T}\int_{0}^{T}e^{-\lambda_{n}\int_{u}^{T}a(s)ds}  e^{-\lambda_{n}\int_{\tau}^{T}a(s)ds}\vert u-\tau \vert^{2H-2} dud\tau \\
&\leq \alpha_{H}\int_{0}^{T}\int_{0}^{T}e^{-\lambda_{n}a_{0}(T-u)}e^{-\lambda_{n}a_{0}(T-\tau)}\vert u-\tau \vert^{2H-2} dud\tau  \\
&\lesssim  \int_{0}^{T-t_{*}} \int_{0}^{T-t_{*}} e^{-\lambda_{n}a_{0}(T-u)}e^{-\lambda_{n}a_{0}(T-\tau)}\vert u-\tau \vert^{2H-2} dud\tau \\
&\quad +\int_{T-t_{*}}^{T} \int_{T-t_{*}}^{T}e^{-\lambda_{n}a_{0}(T-u)}e^{-\lambda_{n}a_{0}(T-\tau)}\vert u-\tau \vert^{2H-2} dud\tau  \\
&\quad +\int_{T-t_{*}}^{T} \int_{0}^{T-t_{*}}e^{-\lambda_{n}a_{0}(T-u)}e^{-\lambda_{n}a_{0}(T-\tau)}\vert u-\tau \vert^{2H-2} dud\tau \\
&\quad +\int_{0}^{T-t_{*}} \int_{T-t_{*}}^{T}e^{-\lambda_{n}a_{0}(T-u)}e^{-\lambda_{n}a_{0}(T-\tau)}\vert u-\tau \vert^{2H-2} dud\tau \\
&:= Q_{1}+Q_{2}+Q_{3}+Q_{4}.
\end{align*}
Substitute variables $u^{'}=T-u$ and $\tau^{'}=T-\tau$ for the above formula we get
\begin{align*}
Q_{1} &= \int_{t_{*}}^{T} \int_{t_{*}}^{T}e^{-\lambda_{n}a_{0}u^{'}}e^{-\lambda_{n}a_{0}\tau^{'}}\vert u^{'}-\tau^{'} \vert^{2H-2} du^{'}d\tau^{'} \\
& \leq \int_{t_{*}}^{T} \int_{t_{*}}^{T}\frac{1}{\lambda_{n}a_{0}u^{'}} \frac{1}{\lambda_{n}a_{0}\tau^{'}}\vert u^{'}-\tau^{'} \vert^{2H-2} du^{'}d\tau^{'} \\
&\lesssim \frac{1}{\lambda_{n}^{2}}\int_{t_{*}}^{T} \int_{t_{*}}^{T}\frac{1}{u^{'}} \frac{1}{\tau^{'}}\vert u^{'}-\tau^{'} \vert^{2H-2} du^{'}d\tau^{'} \\
&\lesssim \frac{1}{\lambda_{n}^{2}}\frac{1}{t_{*}^{2}}\int_{t_{*}}^{T} \int_{t_{*}}^{T} \vert u^{'}-\tau^{'} \vert^{2H-2} du^{'}d\tau^{'} \lesssim \frac{1}{\lambda_{n}^{2}}\frac{1}{t_{*}^{2}},
\end{align*}

\begin{align*}
Q_{2} &= \int_{0}^{t_{*}} \int_{0}^{t_{*}}e^{-\lambda_{n}a_{0}u^{'}}e^{-\lambda_{n}a_{0}\tau^{'}}\vert u^{'}-\tau^{'} \vert^{2H-2} du^{'}d\tau^{'} \\
& \leq \int_{0}^{t_{*}} \int_{0}^{t_{*}}\vert u^{'}-\tau^{'} \vert^{2H-2} du^{'}d\tau^{'} \\
&\lesssim t_{*}^{2H},
\end{align*}

\begin{align*}
Q_{3} &= \int_{0}^{t_{*}} \int_{t_{*}}^{T}e^{-\lambda_{n}a_{0}u^{'}}e^{-\lambda_{n}a_{0}\tau^{'}}\vert u^{'}-\tau \vert^{2H-2} du^{'}d\tau^{'} \\
& \leq \int_{0}^{t_{*}} \int_{t_{*}}^{T}\frac{1}{\lambda_{n}a_{0}u^{'}} \vert u^{'}-\tau^{'} \vert^{2H-2} du^{'}d\tau^{'} \\
&\lesssim \frac{1}{\lambda_{n}t_{*}}\int_{0}^{t_{*}} \int_{t_{*}}^{T}\vert u^{'}-\tau^{'} \vert^{2H-2} du^{'}d\tau^{'} \\
&\lesssim \frac{T^{2H-2}}{\lambda_{n}t_{*}}(T-t_{*})t_{*}\lesssim \frac{1}{\lambda_{n}}.
\end{align*}
Due to the symmetry of $u$ and $\tau$, it is easy to obtain that $Q_{4}=Q_{3}$. From the above estimates about $Q_{1},Q_{2},Q_{3},Q_{4}$ and chose $t_{*}=\lambda_{n}^{-\frac{1}{2}}$ we get that
$$\mathbb{E}\left[ \left( \int_{0}^{T}e^{-\lambda_{n}\int_{\tau}^{T}a(s)ds}d B^{H}(\tau) \right) ^{2}\right] \lesssim  \frac{1}{\lambda_{n}^{2}}\frac{1}{t_{*}^{2}} +t_{*}^{2H}+ \frac{1}{\lambda_{n}}+\frac{1}{\lambda_{n}} \lesssim \frac{1}{\lambda_{n}}.$$

From the above estimates for $H \in (0,1)$, we can conclude that
\begin{align}
\mathbb{E}\left[ \left( \int_{0}^{T}e^{-\lambda_{n}\int_{\tau}^{T}a(s)ds}d B^{H}(\tau) \right) ^{2}\right] \lesssim \frac{1}{\lambda_{n}}+\frac{1}{\lambda_{n}^{2H}} \quad \rightarrow 0\quad as \quad n\rightarrow \infty ,
\end{align}

Combining (4.6) and (4.7) can demonstrate the inverse problem is unstable to recover $f$ and $\vert g\vert$.

\section{Numerical experiments}
\noindent In this section, we present some numerical experiments for the one-dimensional problem where $D = [0,\pi].$ Since the difficulty in getting the exact solution of the problem (1.1), it is assumed that $f(r),$ $h(t),$  and $g(r)$ in the source term are known, and the direct problem is solved by the finite difference method, the final value data $u(r,T)$ can be obtained.  Adding the random disturbance to it, we obtain the noisy data
$$u^{\delta}(r,T) = u(r,T)+ \epsilon u(r,T) \cdot \left(2rand\Big(size\big(u(r,T)\big)\Big)-1\right) $$
where $\epsilon >0$ is relative error level, $\delta = \Vert u^{\delta}(r,T)-u(r,T)  \Vert .$

For some fixed integers $M$ and $N$, we define the time and space step-sizes $h_{t} = T /M$, $h_{r}=\pi / N $ and nodes
$$t_{n} = nh_{t},\quad n\in \{ 0,1,\cdots,M \}, \qquad r_{i}=ih,\quad i\in \{ 0,1,\cdots,N \}.$$

First, we present the numerical solution of the direct problem. We discretize time with the first order backward difference and space with the second order central difference to generate the synthetic data. Let $u_{i}^{n}$ be the numerical approximation to $u(r_{i},t_{n})$. Then we obtain the following explicit scheme:
\begin{align}
\frac{u_{i}^{n}-u_{i}^{n-1}}{h_{t}} &= a(t_{n})\left( \frac{u_{i+1}^{n}-2u_{i}^{n}+u_{i-1}^{n}}{h^{2}}+\frac{2}{r_{i}}\frac{u_{i+1}^{n}-u_{i-1}^{n}}{2h}\right)\nonumber \\
&+f(r_{i})h(t_{n})+g(r_{i})\frac{B^{H}(t_{n})-B^{H}(t_{n-1})}{t_{n}-t_{n-1}}.
\end{align}
where $i=1,2,\cdots,N-1,$ $n=1,2,\cdots,M.$ we have the matrix equation corresponding to the discrete format
$$AU^{n}= \frac{h^{2}}{h_{t}}U^{n-1}+ h^{2}F+\frac{h^{2}}{h_{t}}B+C, \quad n=1,2,\cdots,M,$$
where
\begin{align*}
&F = (f(r_{1})h(t_{n}),f(r_{2})h(t_{n}),\cdots,f(r_{N-1})h(t_{n}))^{T}, \\
&B = \Big(g(r_{1})\left[B^{H}(t_{n})-B^{H}(t_{n-1})\right],g(r_{2})\left[B^{H}(t_{n})-B^{H}(t_{n-1})\right],\cdots,g(r_{N-1})\left[B^{H}(t_{n})-B^{H}(t_{n-1})\right]\Big)^{T},\\
&C =\left(-a(t_{n})\left(\frac{h}{r_{1}}-1\right)u_{0}^{n},0,\cdots,0,a(t_{n})\left(\frac{h}{r_{N-1}}+1\right)u_{N}^{0}\right) \\
&U^{n}= (u_{1}^{n},u_{2}^{n},\cdots,u_{N-1}^{n})^{T}, \quad n=1,2,\cdots,M-1,
\end{align*}
Matrix A is a tridiagonal matrix
\begin{equation}
\left[
\begin{array}{ccccc}
\frac{h^{2}}{h_{t}}+2a(t_{n}) & -a(t_{n})\left(\frac{h}{r_{1}}+1\right) & 0  & 0 & \cdots \\
a(t_{n})\left(\frac{h}{r_{2}}-1\right) & \frac{h^{2}}{h_{t}}+2a(t_{n}) & -a(t_{n})\left(\frac{h}{r_{2}}+1\right) & 0 & \cdots \\
\vdots & \vdots & \vdots  & \ddots & \vdots  \\
\cdots & 0 & a(t_{n})\left(\frac{h}{r_{N-2}}-1\right)  & \frac{h^{2}}{h_{t}}+2a(t_{n}) & -a(t_{n})\left(\frac{h}{r_{N-2}}+1\right) \\
\cdots & 0 & 0 & a(t_{n})\left(\frac{h}{r_{N-1}}-1\right) & \frac{h^{2}}{h_{t}}+2a(t_{n}) \nonumber
\end{array}
\right ]
\end{equation}

Next, we give the numerical solution of the inverse problem. For simplicity, we consider the reconstruction of $g^{2}$
instead of $\left| g\right|$. Recall that functions $f,g \in L^{2}(D)$ have the expanisions
\begin{align}
f=\sum_{n=1}^{\infty}f_{n}\omega_{n},\qquad g^{2}=\sum_{m,n=1}^{\infty}g_{m}g_{n}\omega_{m}\omega_{n},
\end{align}
and coefficients $f_{n}$ and $g_{m}g_{n}$ can be recovered by (4.1) and (4.3), respectively,
\begin{align}
f_{n}&=\frac{\mathbb{E}(u_{n}(T))}{\int_{0}^{T}h(\tau)e^{-\lambda_{n}\int_{\tau}^{T}a(s)ds}d\tau}\nonumber, \\
g_{m}g_{n}&=\frac{Cov(u_{m}(T),u_{n}(T))}{\mathbb{E}\left[ \int_{0}^{T}e^{-\lambda_{m}\int_{\tau}^{T}a(s)ds}d B^{H}(\tau)\int_{0}^{T}e^{-\lambda_{n}\int_{\tau}^{T}a(s)ds}d B^{H}(\tau) \right]},
\end{align}
where apply composite Simpson rule to calculate the integral. Since the inverse problem is ill-posed, we truncate above series by the first $N_{1}$ terms as a regularization. We use the parameter P to denote the number of sample paths, which is used to approximate the expectation involved in recovery formula (5.3).

In the numerical experiments, we choose $N_{1}=30,$ $M=2^{11},$ $N=100,$ and $T=1,$ where the temporal step is meant to be small such that the numerical solution is accurate enough for the direct problem. Functions in (1.1) are chosen as $a(t)=t^2,$ $f(r)=\sin (3r),$ $h(t)=1.$ The total number of sample paths $P=1000$ are used when simulating the covariance of the solution. In addition, the data is assumed to be polluted by a uniformly distributed noise with the error level $\epsilon = 0.001$.\\
{\bf Example 1.} Consider a smooth function
$$g(r)=\sin(2r).$$
The reconstruction results of the $f$ and $g^{2}$ with $H = 0.1$,  $H = 0.5$ and $ H = 0.9$ are given in figures $1-3$, respectively.

\begin{figure}
  \centering
  \includegraphics[width=0.4\textwidth]{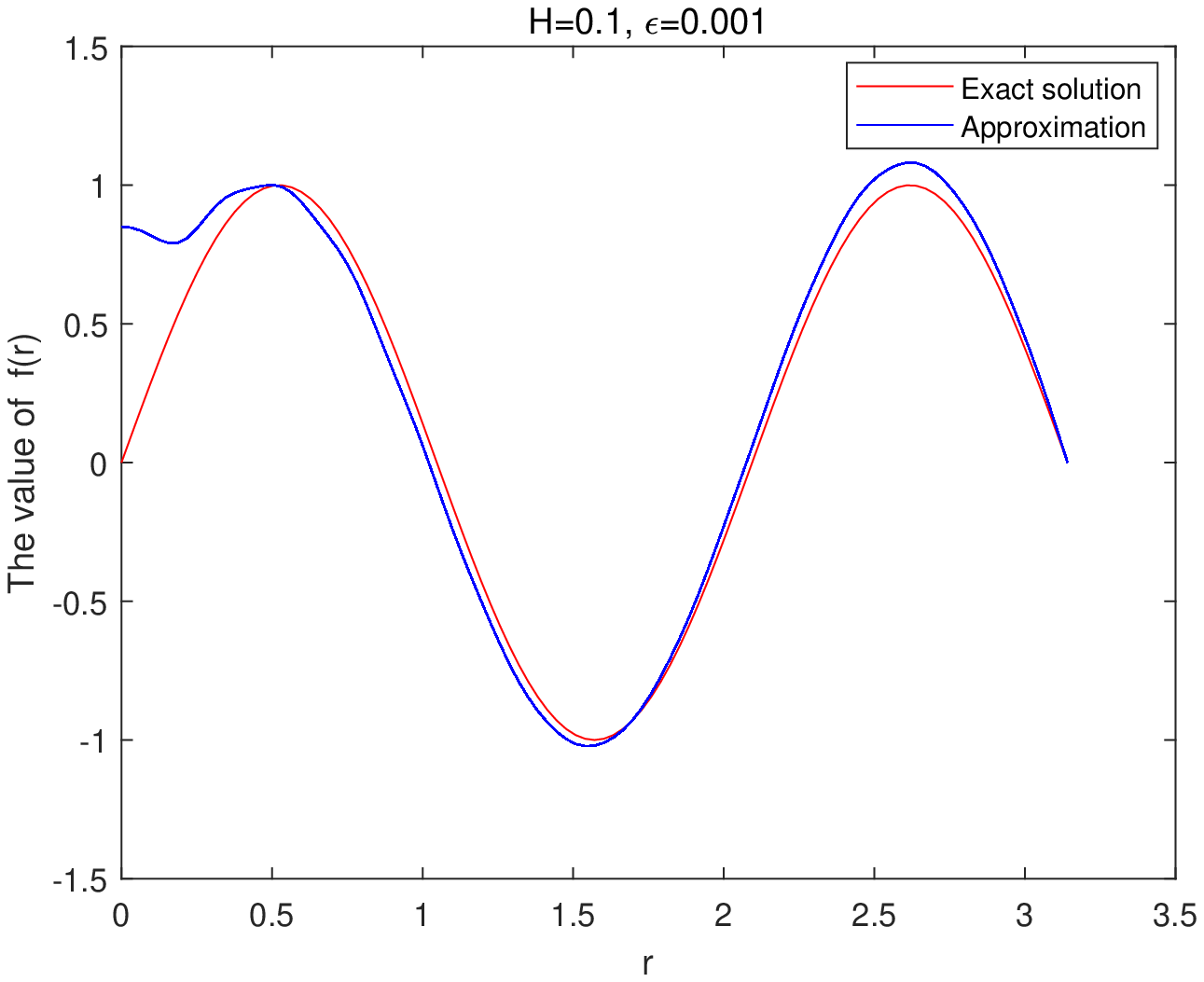}
  \includegraphics[width=0.4\textwidth]{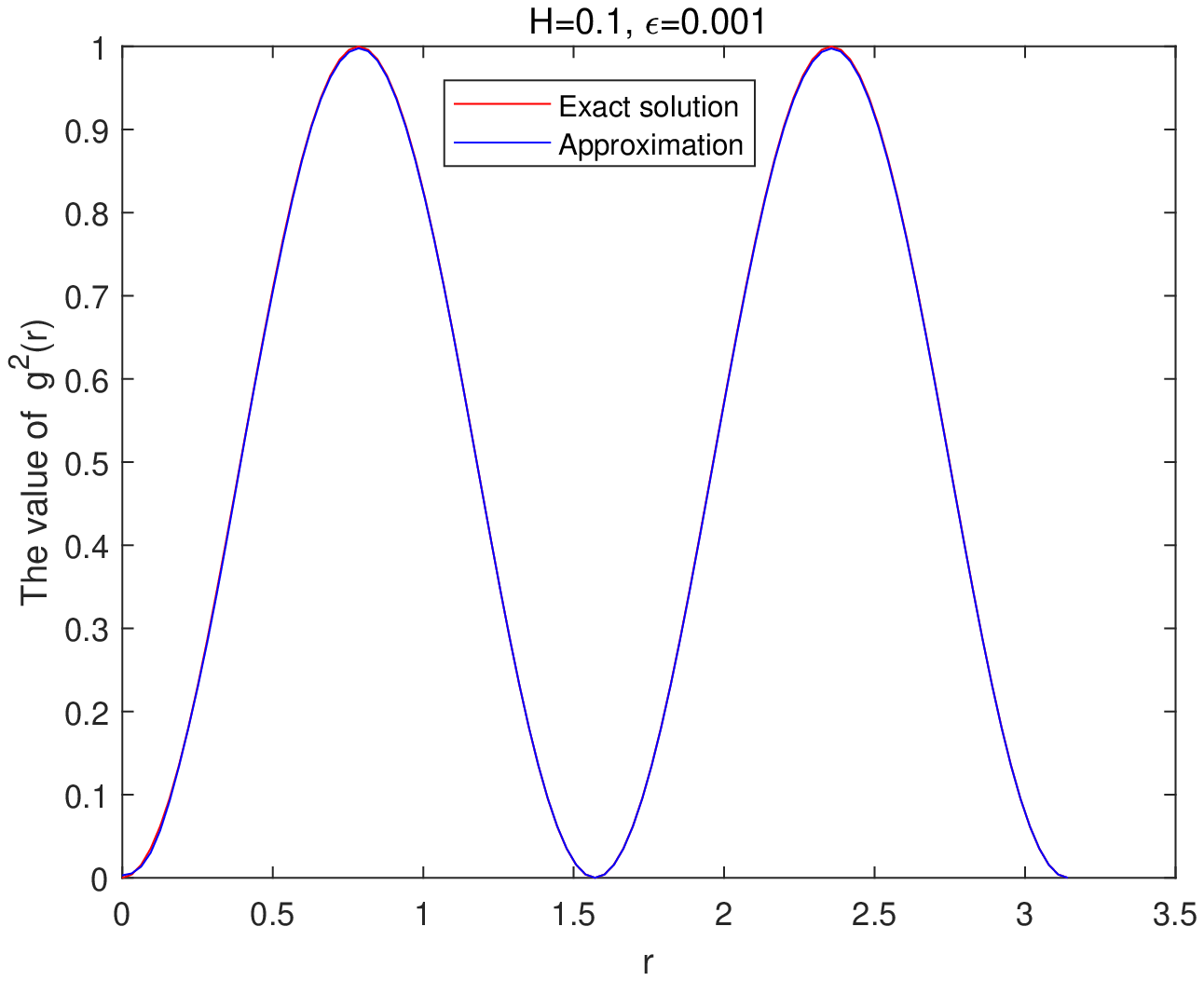}
\caption{The reconstruction for $f$ (left column) and $g^{2}$ (right column) for $H=0.1$ with $\epsilon=0.001.$}
\end{figure}

\begin{figure}
  \centering
  \includegraphics[width=0.4\textwidth]{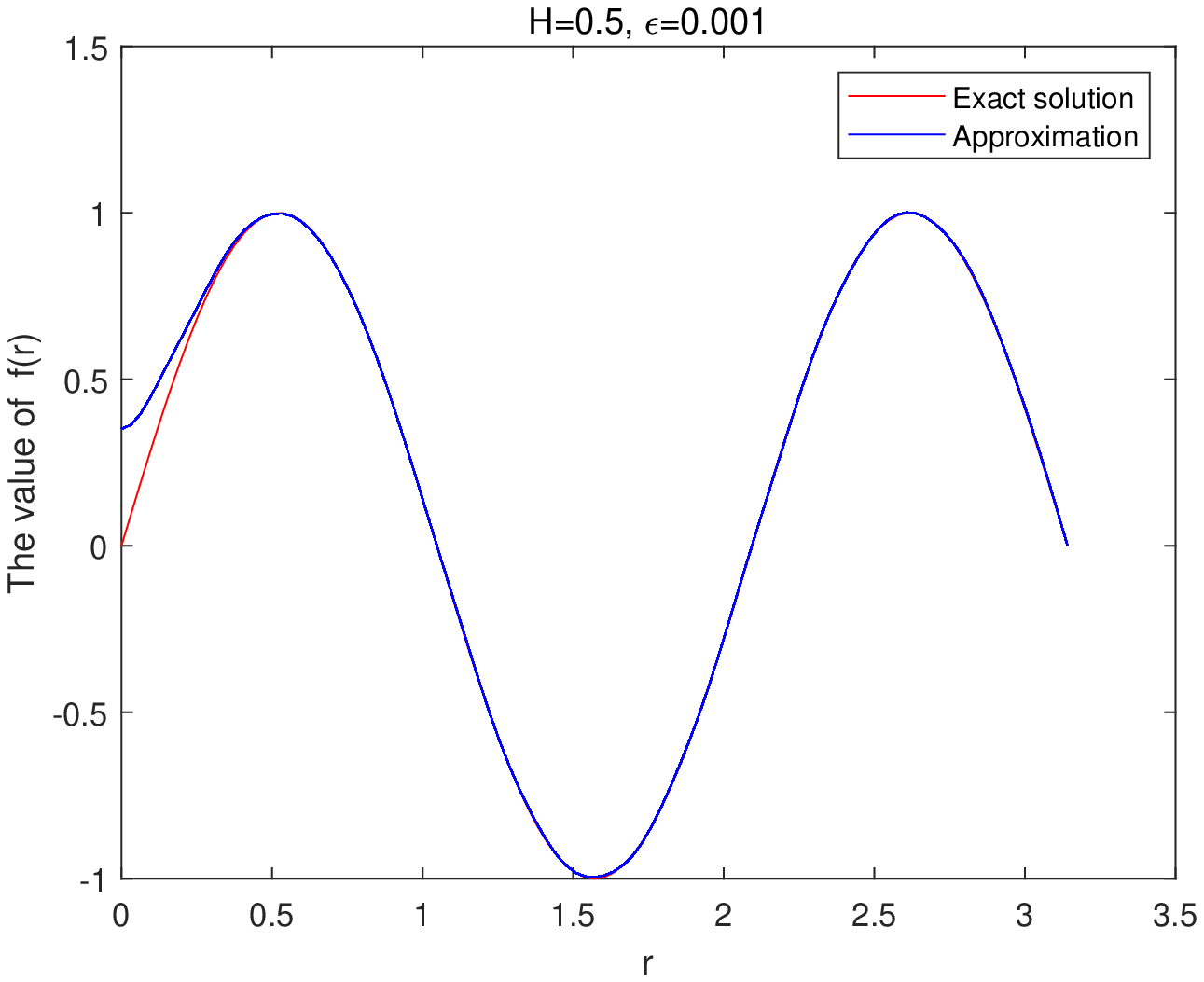}
  \includegraphics[width=0.4\textwidth]{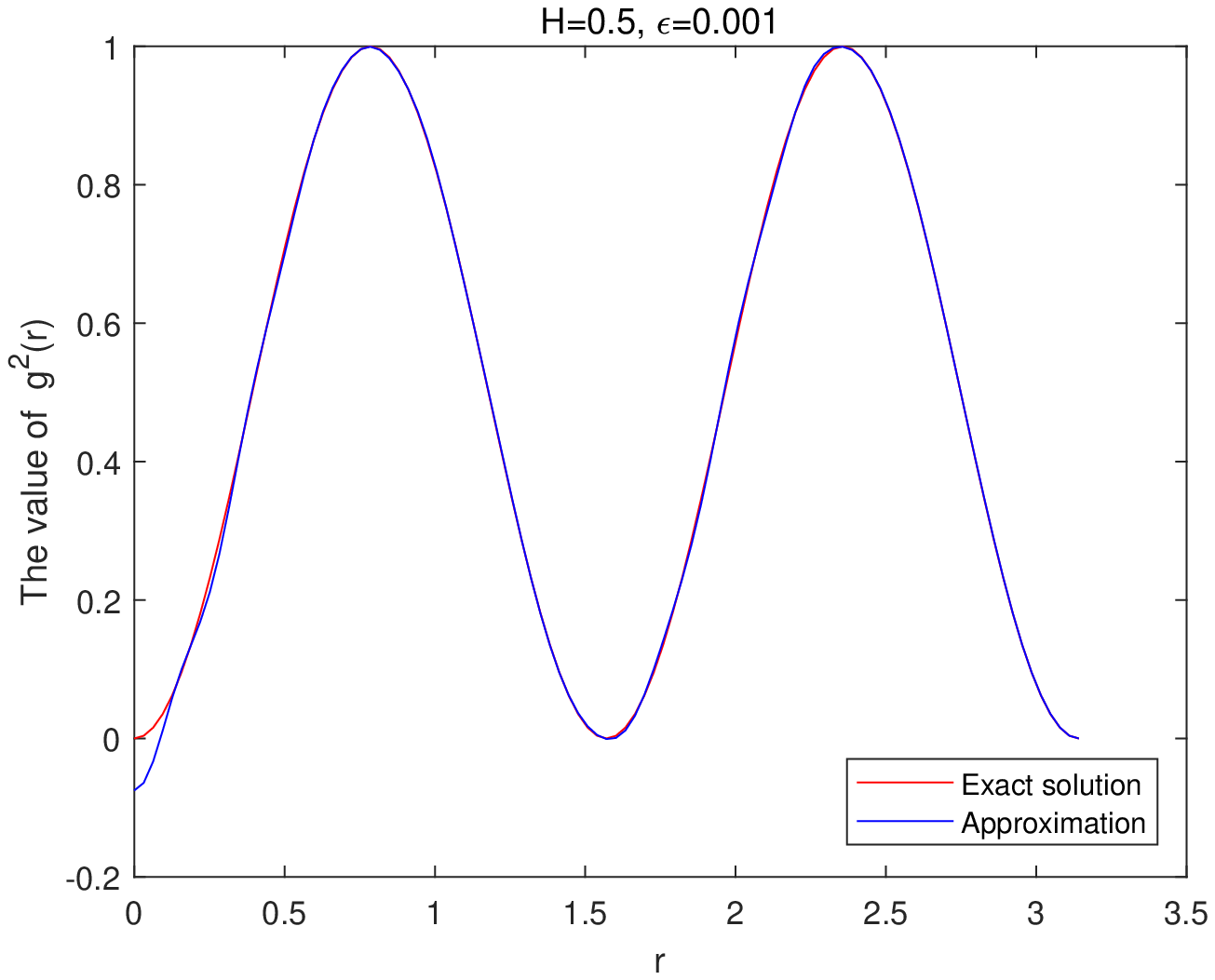}
\caption{The reconstruction for $f$ (left column) and $g^{2}$ (right column) for $H=0.5$ with $\epsilon=0.001.$}
\end{figure}

\begin{figure}
  \centering
  \includegraphics[width=0.4\textwidth]{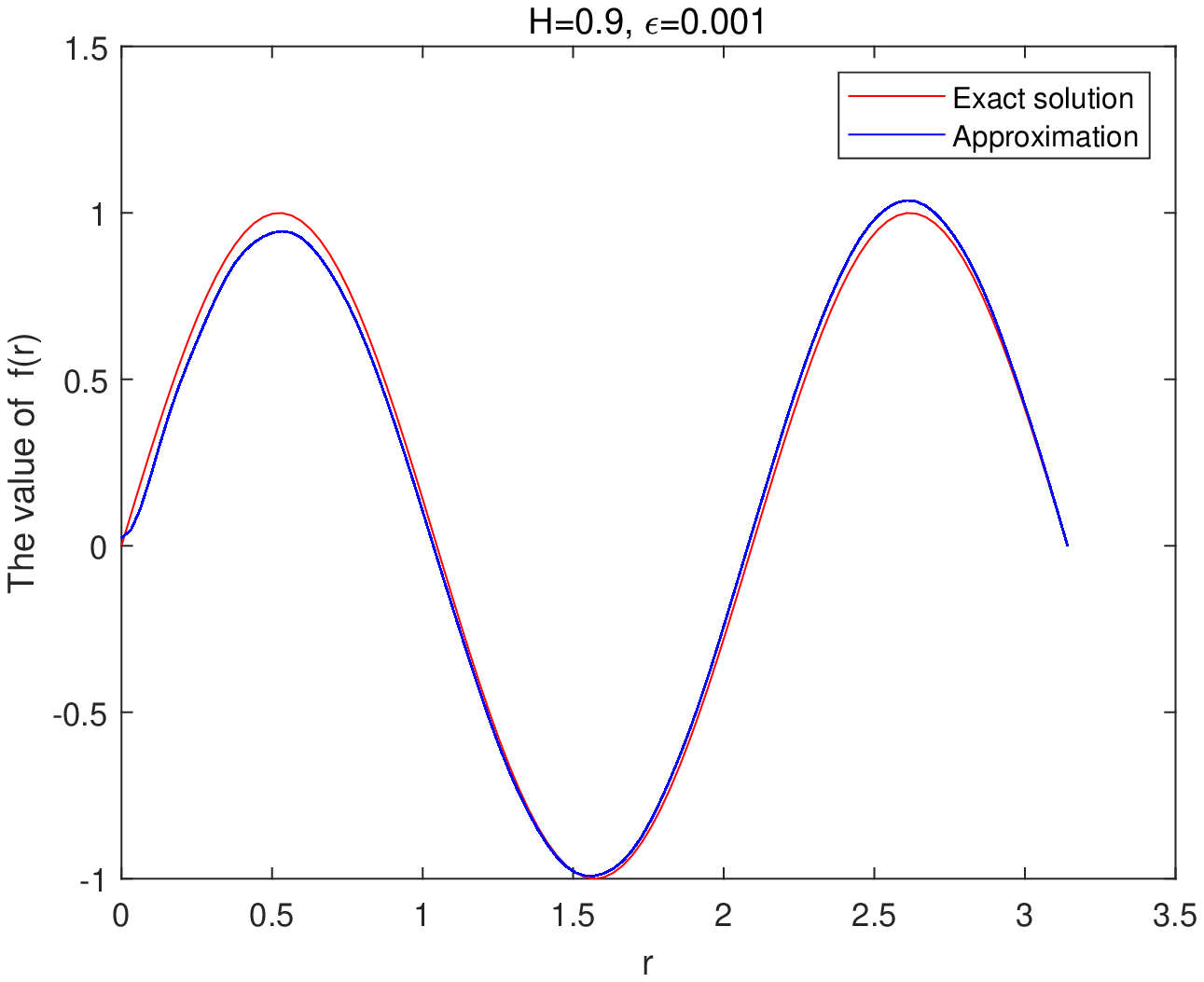}
  \includegraphics[width=0.4\textwidth]{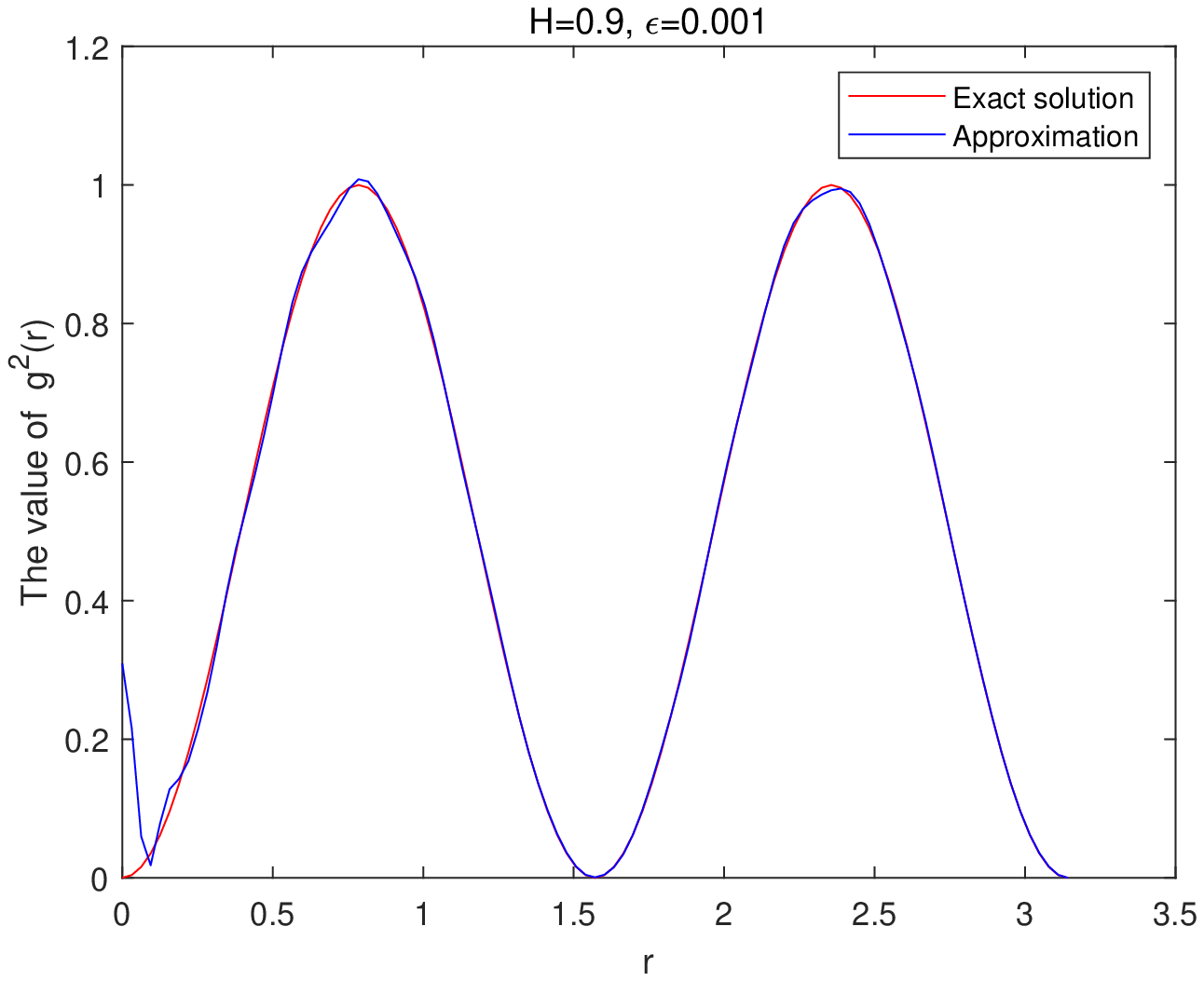}
\caption{The reconstruction for $f$ (left column) and $g^{2}$ (right column) for $H=0.9$ with $\epsilon=0.001.$}
\end{figure}

\noindent {\bf Example 2.} In this case, we consider a piecewise smooth function
$$g(r)=\left\{
\begin{aligned}
&0,                  &0\leq r\leq \frac{1}{4}\pi ,\\
&4r-\pi,\qquad &\frac{1}{4}\pi < r\leq \frac{1}{2}\pi ,\\
&3\pi-4r,\qquad&\frac{1}{2}\pi < r\leq \frac{3}{4}\pi ,\\
&0,                 &\frac{3}{4}\pi< r\leq \pi.
\end{aligned}
\right.
$$
For this case, reconstructions of the $f$ and $g^{2}$ are showed in Figure $4-6$, in which the Hurst index $H$ is taken as 0.1 ,0.5 and 0.9 respectively.

\begin{figure}
  \centering
  \includegraphics[width=0.4\textwidth]{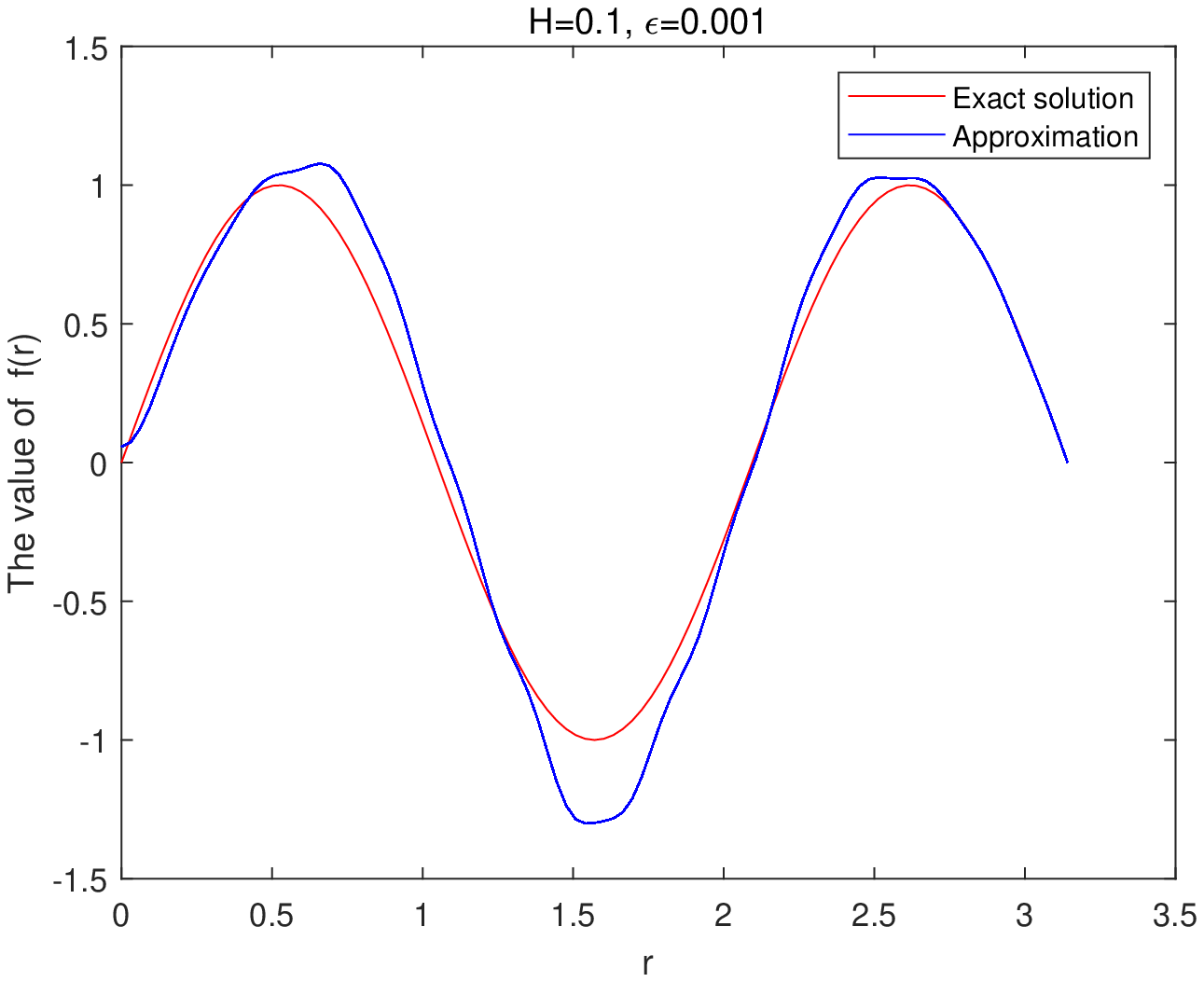}
  \includegraphics[width=0.4\textwidth]{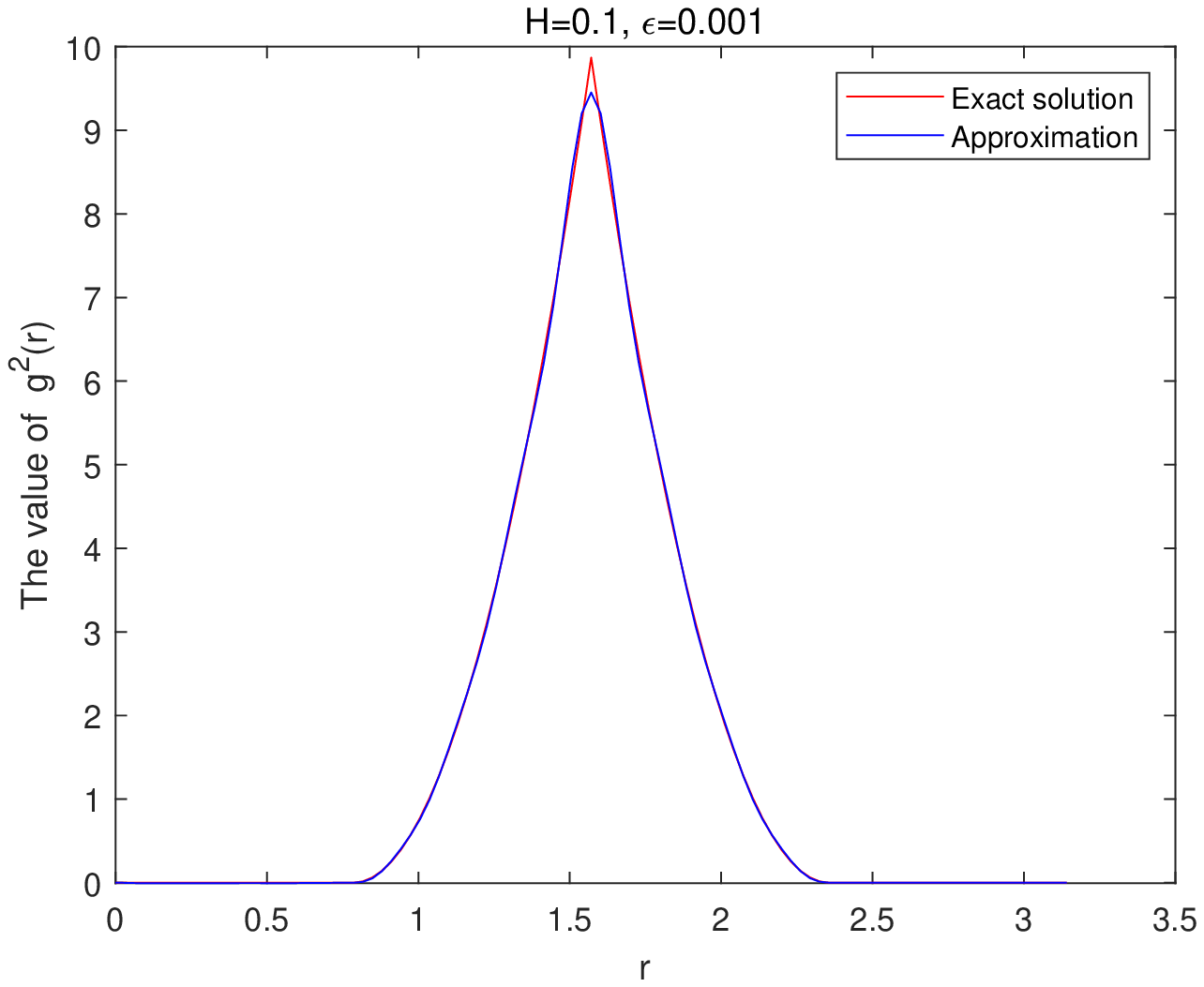}
\caption{The reconstruction for $f$ (left column) and $g^{2}$ (right column) for $H=0.4$ with $\epsilon=0.001.$}
\end{figure}

\begin{figure}
   \centering
  \includegraphics[width=0.4\textwidth]{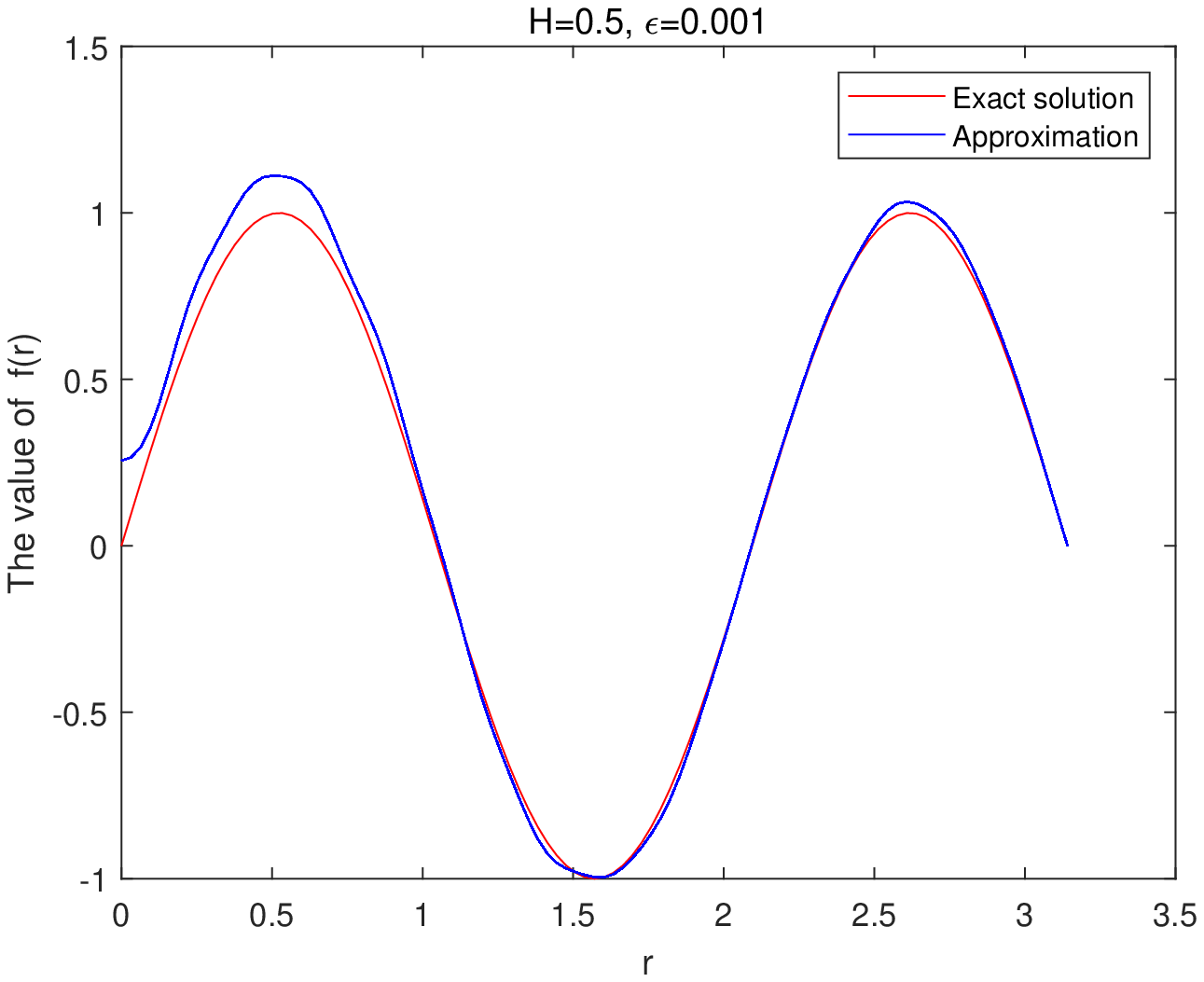}
  \includegraphics[width=0.4\textwidth]{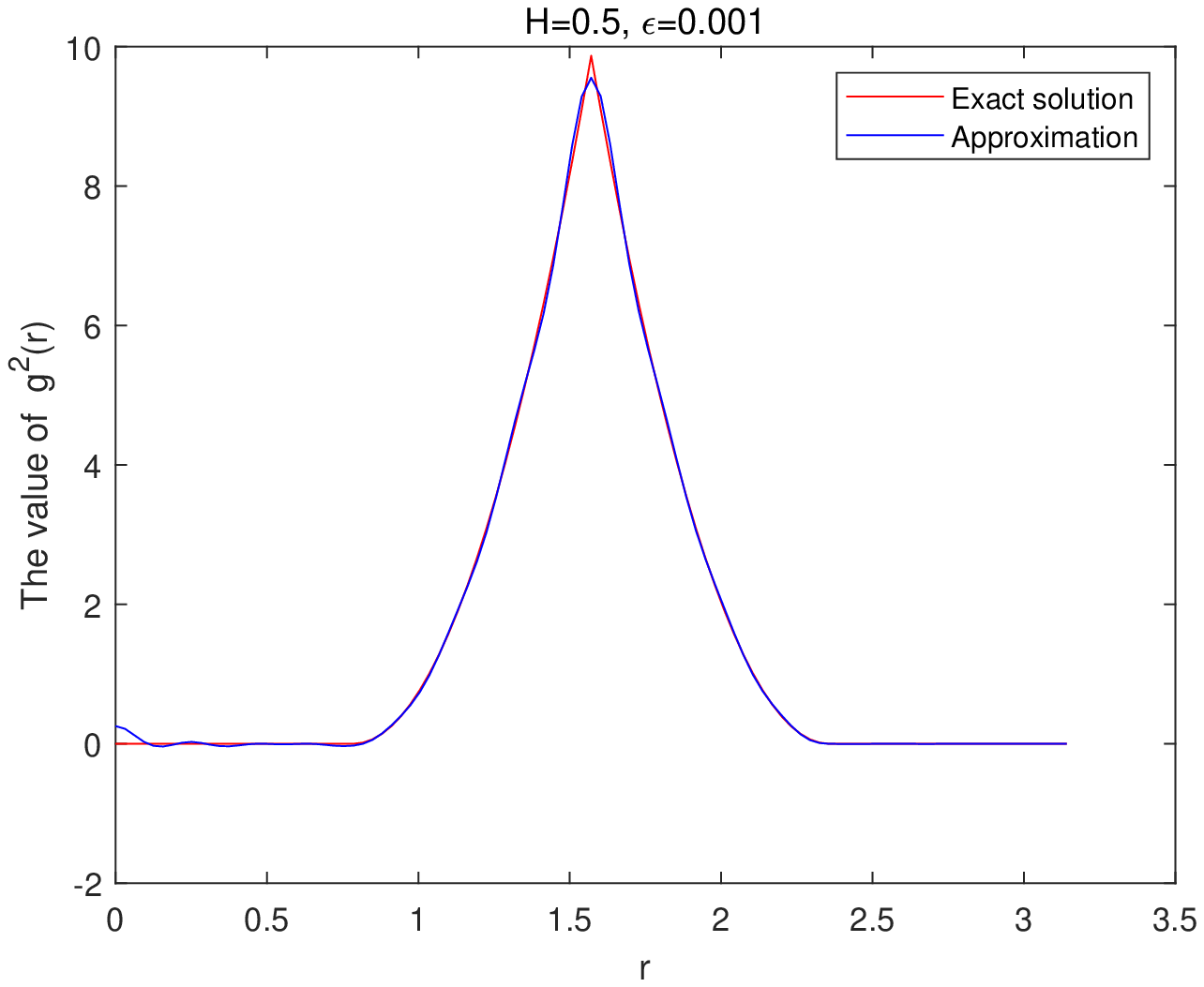}
\caption{The reconstruction for $f$ (left column) and $g^{2}$ (right column) for $H=0.5$ with $\epsilon=0.001.$}
\end{figure}

\begin{figure}
  \centering
  \includegraphics[width=0.4\textwidth]{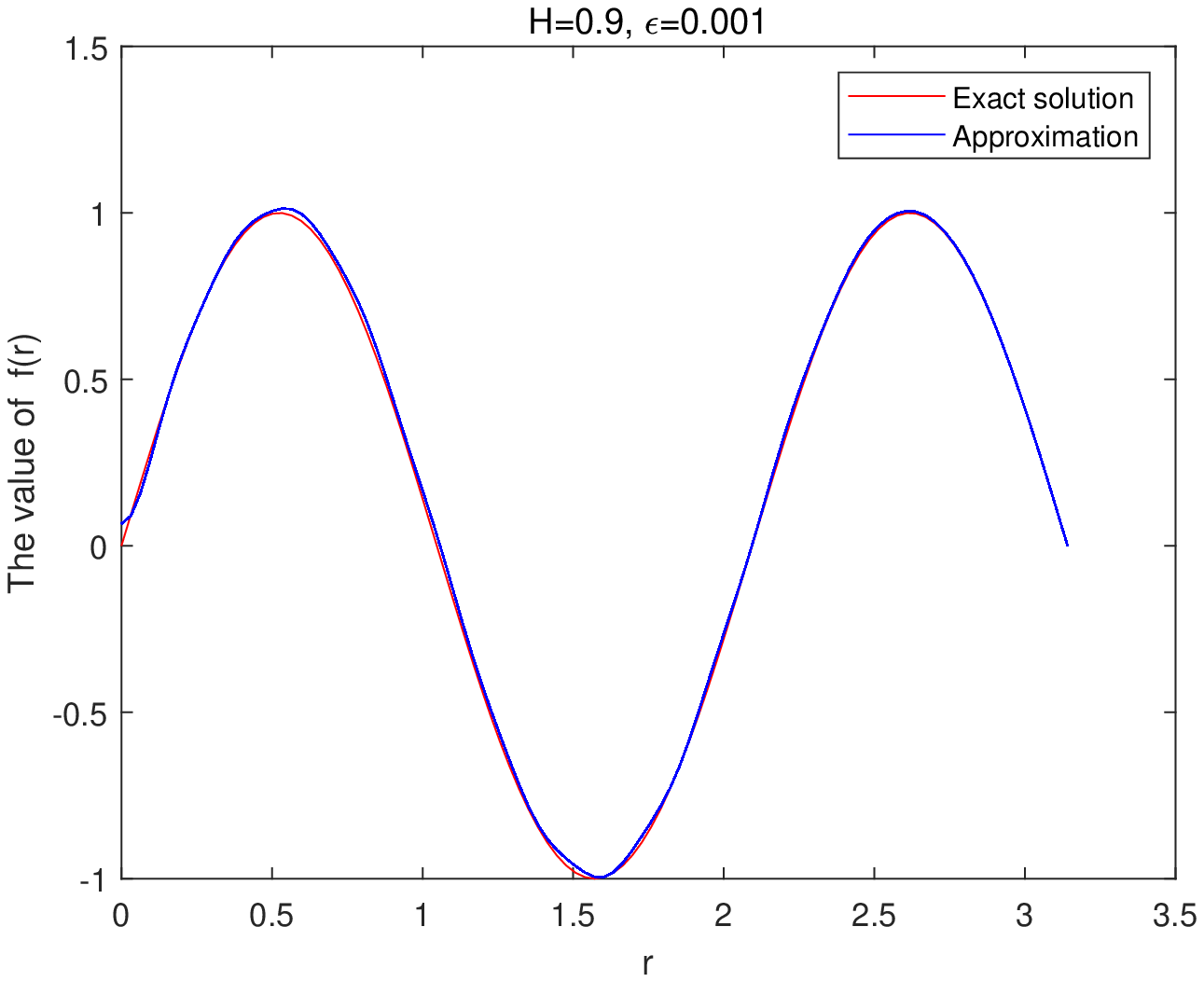}
  \includegraphics[width=0.4\textwidth]{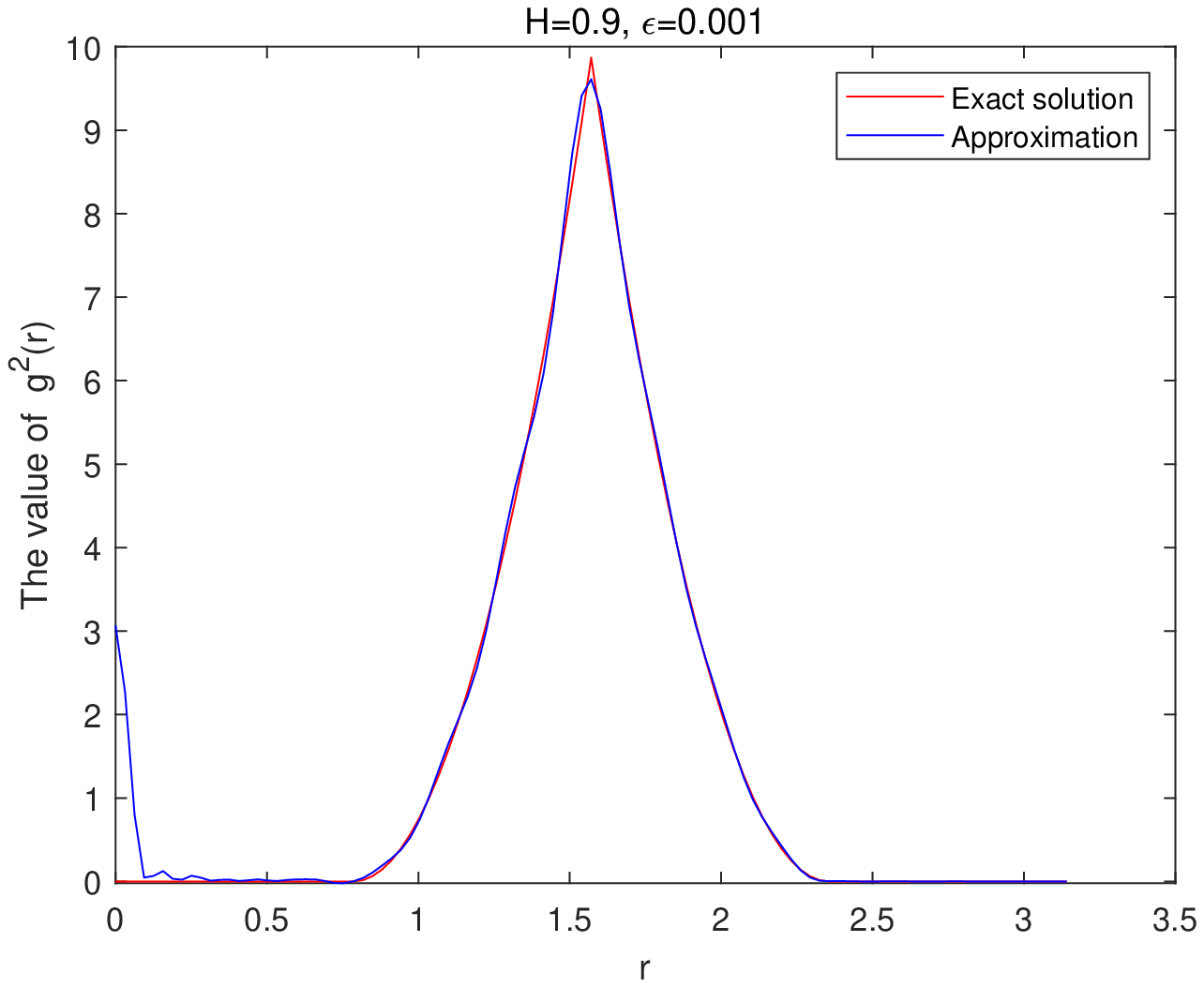}
\caption{The reconstruction for $f$ (left column) and $g^{2}$ (right column) for $H=0.9$ with $\epsilon=0.001.$}
\end{figure}

\noindent {\bf Example 3.} Consider a non-smooth piecewise constant function
$$g(r)=\left\{
\begin{aligned}
&0,                 &0\leq r\leq \frac{1}{3}\pi ,\\
&1,          \qquad     &\frac{1}{3}\pi < r\leq \frac{2}{3}\pi ,\\
&0,                 &\frac{2}{3}\pi< r\leq \pi.
\end{aligned}
\right.
$$
The reconstructing results are displayed in Figure $7-9$. Similarly, we choose the Hurst index $H$ as 0.1, 0.5 and 0.9 respectively.

\begin{figure}
  \centering
  \includegraphics[width=0.4\textwidth]{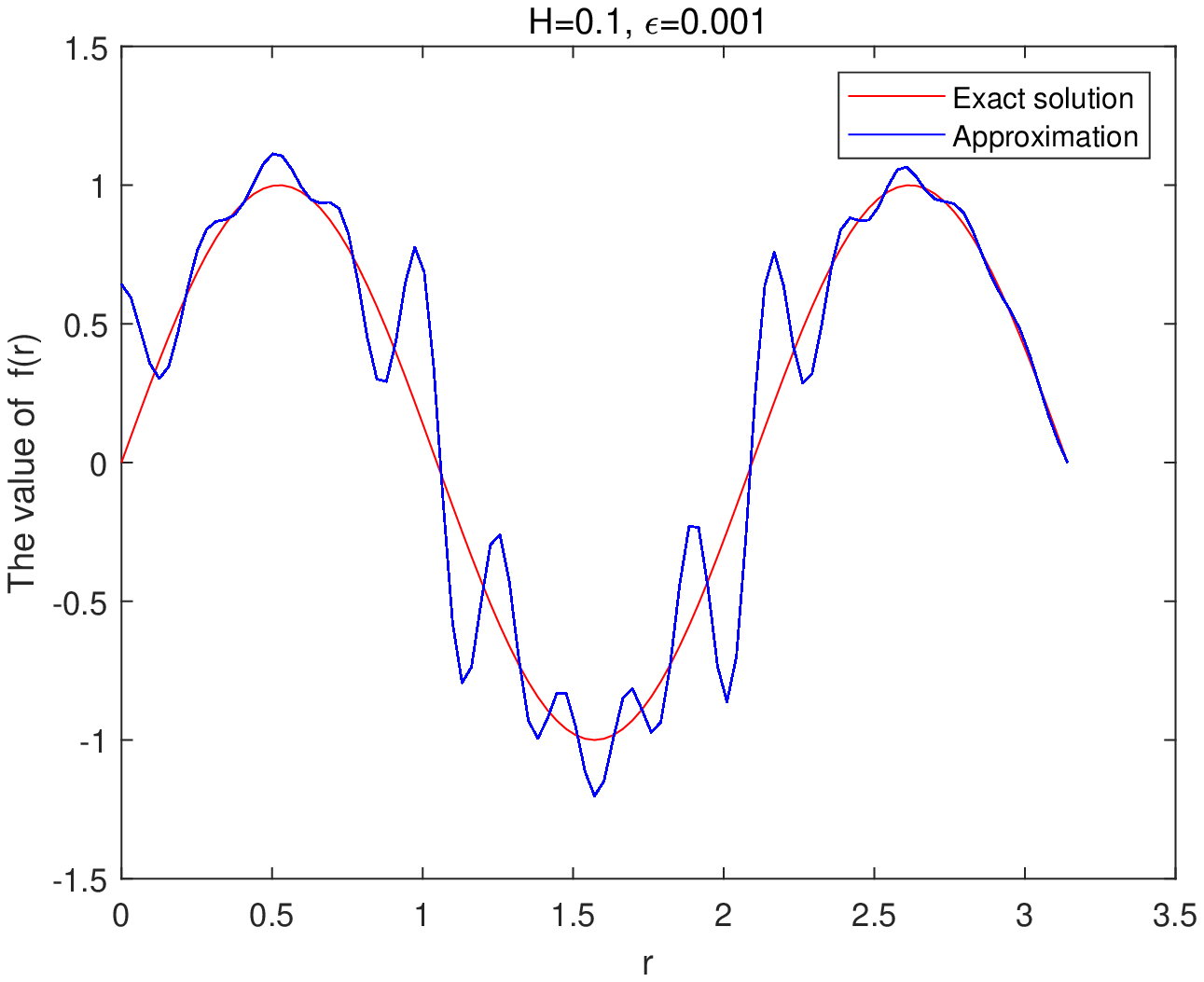}
  \includegraphics[width=0.4\textwidth]{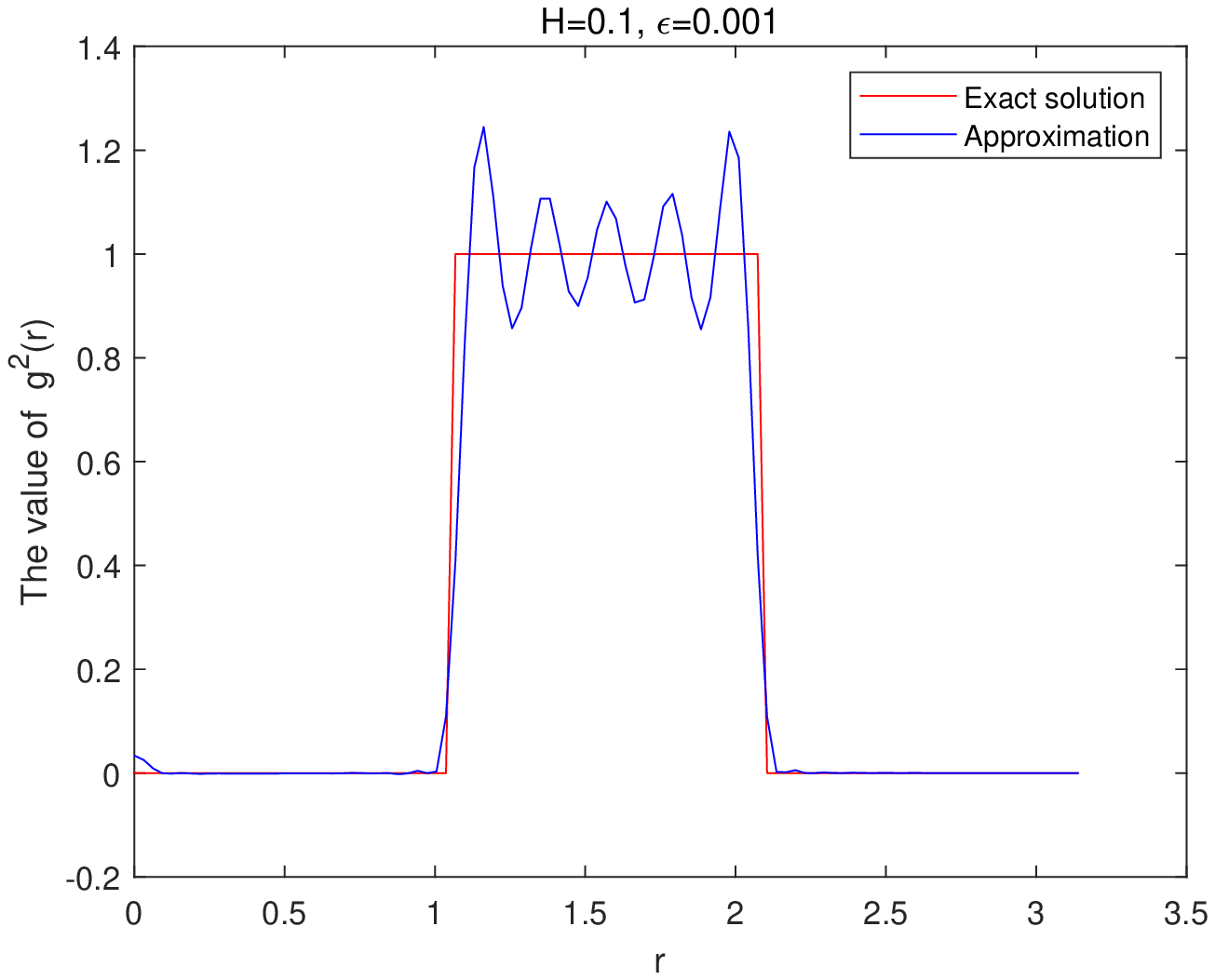}
\caption{The reconstruction for $f$ (left column) and $g^{2}$ (right column) for $H=0.4$ with $\epsilon=0.001.$}
\end{figure}

\begin{figure}
  \centering
  \includegraphics[width=0.4\textwidth]{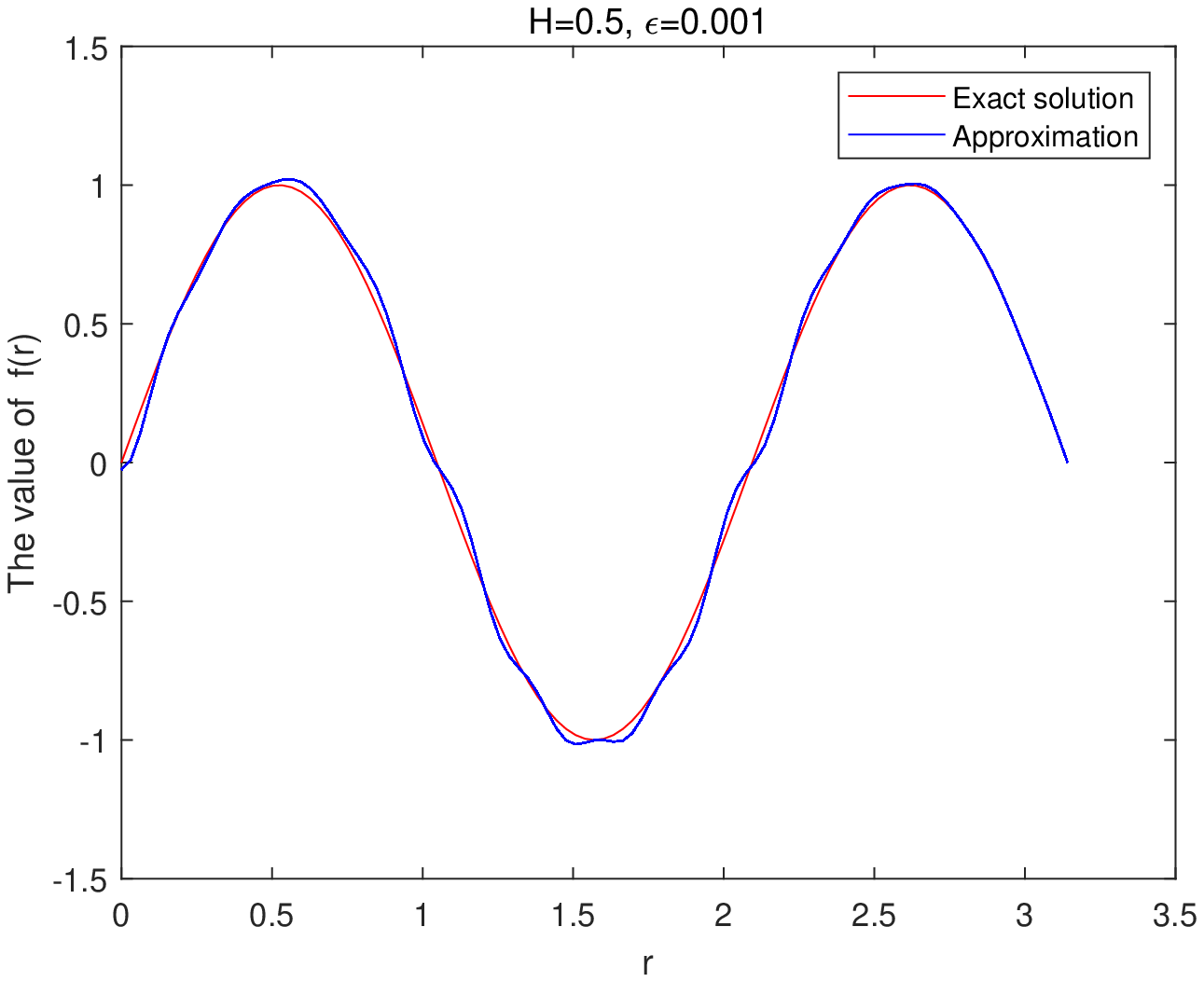}
  \includegraphics[width=0.4\textwidth]{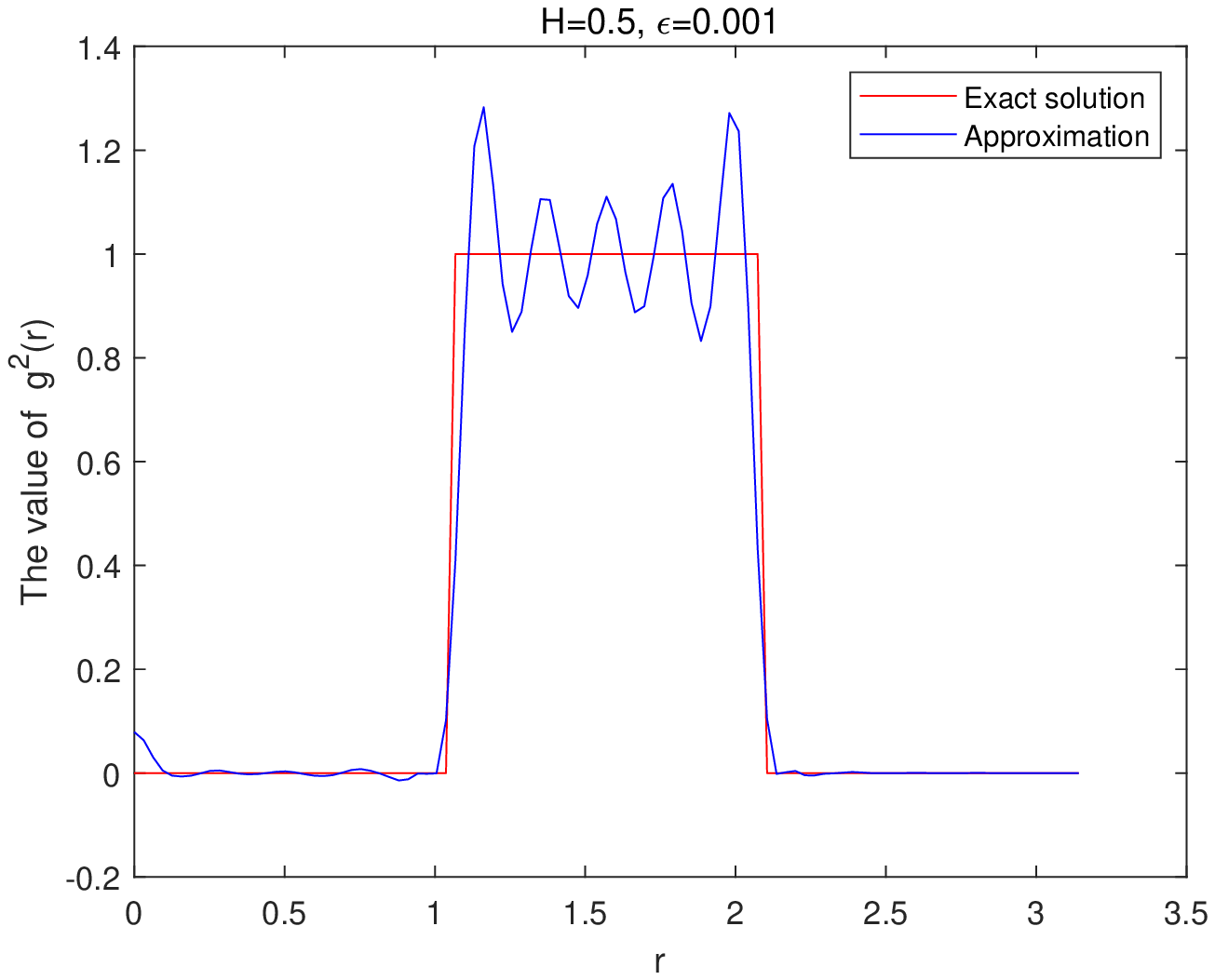}
\caption{The reconstruction for $f$ (left column) and $g^{2}$ (right column) for $H=0.5$ with $\epsilon=0.001.$}
\end{figure}

\begin{figure}
  \centering
  \includegraphics[width=0.4\textwidth]{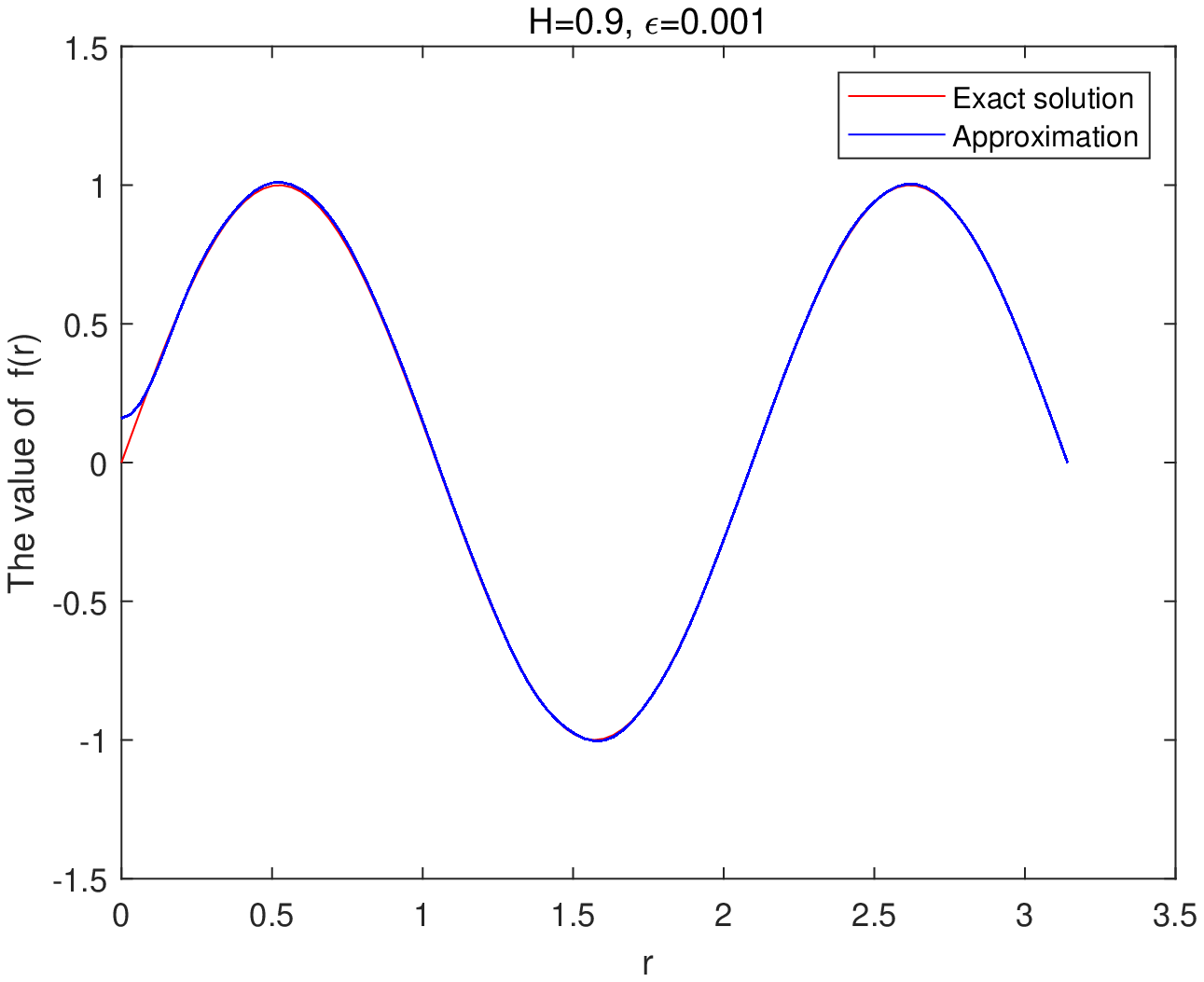}
  \includegraphics[width=0.4\textwidth]{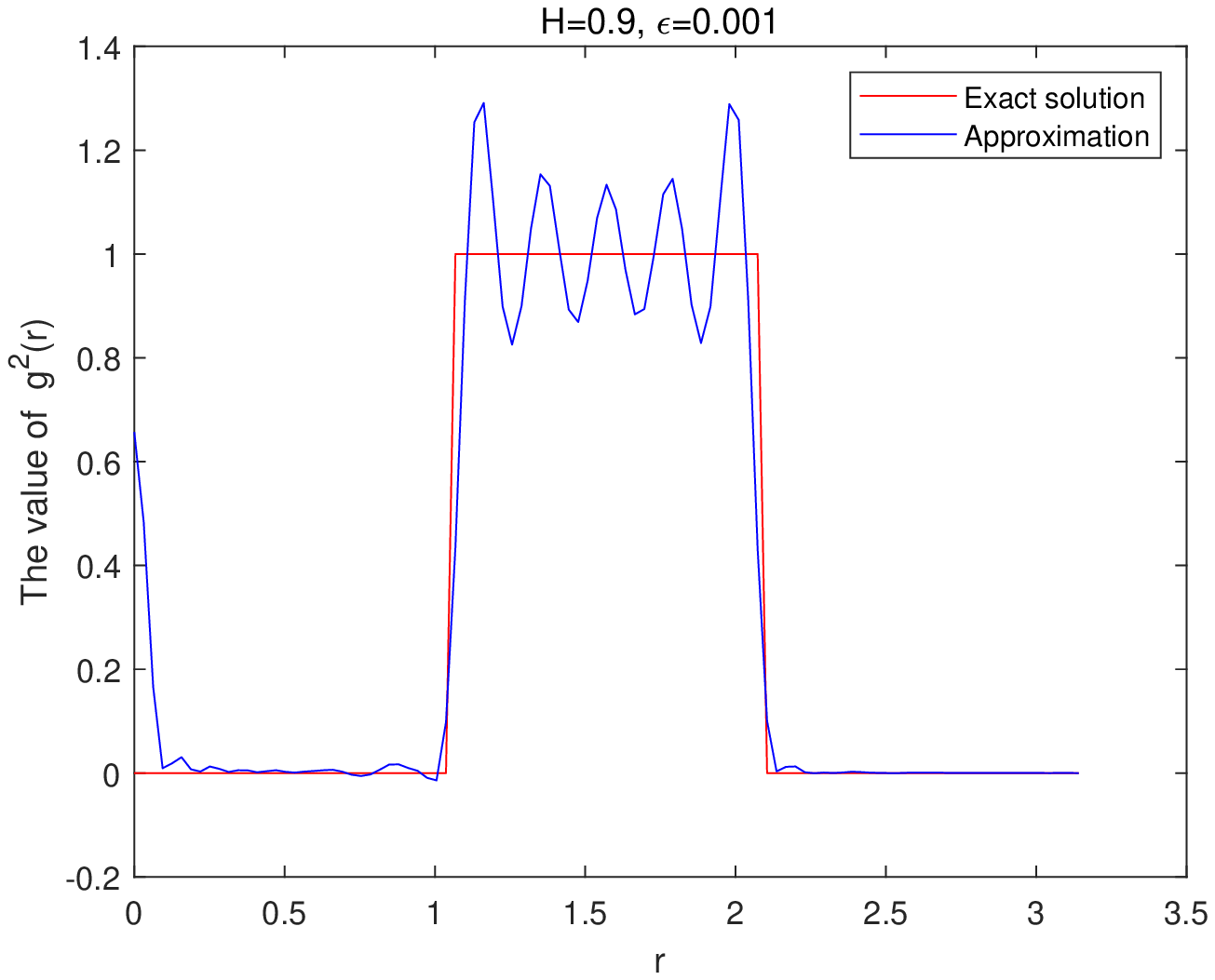}
\caption{The reconstruction for $f$ (left column) and $g^{2}$ (right column) for $H=0.9$ with $\epsilon=0.001.$}
\end{figure}

From Figures $1-9$, we can find that the numerical results of smooth source term are the best, the numerical results of the non-smooth but continuous
source term are slightly worse, and the numerical results of the discontinuous source term are the worst. What is more, comparing figures with different $H$, we can see that the reconstruction would be better if $H$ is larger.

\section{Conclusion}
\noindent In this paper, we study the inverse random source problem of Helium production-diffusion equation driven by fractional Brownian motions. For the direct problem, we obtain some well-posed results. For the inverse source problem, we give the uniqueness and the instability results. In the numerical experiments, three different examples shows the effectiveness and feasibility of the method.\\

\noindent\textbf{Acknowledgments}\\
The authors would like to thank Prof. Xiaoli Feng (School of Mathematics and Statistics, Xidian University)'s help.\\

\end{document}